\documentclass[12pt,a4paper]{amsart}
\usepackage{latexsym, amsfonts, amsmath, amssymb}
\usepackage{fullpage}
\newcommand{\beq}{\begin{equation}}
\newcommand{\eeq}{\end{equation}}
\newcommand{\bea}{\begin{eqnarray}}
\newcommand{\eea}{\end{eqnarray}}
\newcommand{\bean}{\begin{eqnarray*}}
\newcommand{\eean}{\end{eqnarray*}}

\newcommand{\brray}{\begin{array}}
\newcommand{\erray}{\end{array}}
\newcommand{\ben}{\begin{equation}{nonumber}}
\newcommand{\een}{\end{equation}{nonumber}}

\newtheorem{dfn}{Definition}[section]
\newtheorem{thm}[dfn]{Theorem}
\newtheorem{lmma}[dfn]{Lemma}
\newtheorem{ppsn}[dfn]{Proposition}
\newtheorem{crlre}[dfn]{Corollary}
\newtheorem{xmpl}[dfn]{Example}
\newtheorem{rmrk}[dfn]{Remark}
\newtheorem{nota}[dfn]{Notation}
\newtheorem{assu}[dfn]{Assumption}
\newcommand{\bdfn}{\begin{dfn}}
\newcommand{\bthm}{\begin{thm}}
\newcommand{\blmma}{\begin{lmma}}
\newcommand{\bppsn}{\begin{ppsn}}
\newcommand{\bcrlre}{\begin{crlre}}
\newcommand{\bxmpl}{\begin{xmpl}}
\newcommand{\brmrk}{\begin{rmrk}}
\newcommand{\edfn}{\end{dfn}}
\newcommand{\ethm}{\end{thm}}
\newcommand{\elmma}{\end{lmma}}
\newcommand{\eppsn}{\end{ppsn}}
\newcommand{\ecrlre}{\end{crlre}}
\newcommand{\exmpl}{\end{xmpl}}
\newcommand{\ermrk}{\end{rmrk}}

\newcommand{\IC}{\mathbb{C}}

\newcommand{\oneform}{{\Omega}^1 ( \mathcal{A} )}

\newcommand{\twoform}{{\Omega}^2( \mathcal{A} )}
\newcommand{\tensora}{\otimes_{\mathcal{A}}}

\newcommand{\tensorc}{\otimes_{\mathbb{C}}}
\newcommand{\A}{\mathcal{A}}

\newcommand{\E}{\mathcal{E}}
\newcommand{\F}{\mathcal{F}}
\newcommand{\Z}{\mathcal{Z}}
\newcommand{\Acenter}{\mathcal{Z}( \mathcal{A} )}
\newcommand{\Ecenter}{\mathcal{Z}( \mathcal{E} )}

\newcommand{\Psym}{P_{\rm sym}}
\newcommand{\Hom}{{\rm Hom}}
\newcommand{\omegazero}{\omega_{(0)}}
\newcommand{\omegaone}{\omega_{(1)}}
\newcommand{\etazero}{\eta_{(0)}}
\newcommand{\etaone}{\eta_{(1)}}
\newcommand{\zeroomega}{{}_{(0)}\omega}
\newcommand{\oneomega}{{}_{(1)}\omega}
\newcommand{\zeroeta}{{}_{(0)}\eta}
\newcommand{\oneeta}{{}_{(1)}\eta}
\newcommand{\zerotheta}{{}_{(0)}\theta}
\newcommand{\onetheta}{{}_{(1)}\theta}

\newcommand{\ancenter}{\mathcal{Z}(A_N)}
\newcommand{\tensoran}{\otimes_{A_N}}

\def\a*{{\cal A}_{h,*}}
\def\B{{\cal B}(h)}
\def\B1{{\cal B}_1(h)}
\def\b{{\cal B}^{\rm s.a.}(h)}
\def\b1{{\cal B}^{\rm s.a.}_1(h)}

\newcommand{\nn}{\nonumber}

\newcommand{\id}{\mbox{id}}

\def \qed {$\Box$}

\def\a*{{\cal A}_{h,*}}
\def\B{{\cal B}(h)}
\def\B1{{\cal B}_1(h)}
\def\b{{\cal B}^{\rm s.a.}(h)}
\def\b1{{\cal B}^{\rm s.a.}_1(h)}

\newcommand{\thetazero}{\theta_{(0)}}
\newcommand{\thetaone}{\theta_{(1)}}

\newcommand{\RNum}[1]{\uppercase\expandafter{\romannumeral #1\relax}}

\begin{document}
\begin{center}
{\Large{\bf On the Koszul formula in noncommutative geometry}}\\
\vspace{0.5in}
{\large Jyotishman Bhowmick, Debashish Goswami, Giovanni Landi}\\
\end{center}
\date{May 2020}
\title{}
\maketitle
\begin{abstract}
We prove a Koszul formula for the Levi-Civita connection for any pseudo-Riemannian bilinear metric on a class of centered bimodule of noncommutative one-forms. As an application to the Koszul formula, we show that our Levi-Civita connection is a bimodule connection. We construct a spectral triple on a fuzzy sphere and compute the scalar curvature for the Levi-Civita connection associated to a canonical metric.
\end{abstract}

\vspace{0.5in}
\tableofcontents
\parskip = 1 ex

\section{Introduction}
The relatively recent paper \cite{Connes_moscovici}  (following earlier work in \cite{scalar_3}) on the computation of 
the scalar curvature for a class of metrics on noncommutative tori, and related Gauss--Bonnet theorems, has lead to a  
flourishing of what one might call noncommutative Riemannian geometry. An approach to this, for which we  refer to the expository paper \cite{khalkhali} and references therein, is via the spectral properties of the Laplacian (or Dirac) operator. 

On the other hand, an algebraic approach was taken by a number of authors (\cite{rosenberg,sheu,article1,article2,pseudo,cylinder,tiger}) whereas one computes the curvature in terms of the Levi-Civita connection. A somewhat related approach is in \cite{majid_1,majid_2,majidkoszul2}. Some earlier studies of metric approach and Levi-Civita connection for noncommutative spaces are in \cite{LaMa88,frolich,LNW94,dubois,dubois2,DHLS96,heckenberger}.

 The proof of the uniqueness of Levi-Civita connections in classical differential geometry yields the Koszul formula and it turns out that this formula actually defines a connection which is torsionless and compatible with the metric. The goal of this article is to demonstrate a noncommutative analogue of this proof under some reasonable assumptions (see Theorem \ref{koszul21stjuly} and Theorem \ref{existenceuniqueness}).
   
Connections in noncommutative geometry have been studied from several viewpoints. In \cite{dubois,dubois2} there were studied covariant derivatives on a certain class of modules of derivations of a noncommutative algebra as well as the notion of bimodule connections (see for these also \cite{beggsmajidbook} and references therein).   Now in classical differential geometry, the Riemannian metric and connections live more naturally  on the level of vector fields. However, in the context of noncommutative geometry, it seems more natural to work on the level of differential forms and that is what we do in the present paper. Thus, for a (possibly) noncommutative algebra $\A$ and the bimodule of one-forms $\E$ coming from a differential calculus, a (right) connection on $ \E $ will be a $\mathbb{C}$-linear map $ \nabla: \E \rightarrow \E \tensora \E $ satisfying the Leibniz rule for the right multiplication of elements in $\A$. However, our approach allows us to make contact with both of these approaches. In Section \ref{section7}, we prove that the Levi-Civita connection that we obtain in Theorem \ref{existenceuniqueness} is a bimodule connection with respect to a canonical symmetrization map obtained from natural assumptions. In a companion article \cite{article4}, we prove that our assumptions allow us to have a sufficiently large $ \Acenter $-bimodule of derivations on the $\A$-bimodule $\E$ of forms so that we can define covariant derivatives and recover a Koszul formula on this $\Acenter$-bimodule.

Let us discuss the plan of the article. In Section \ref{13thjuly20191}, we recall the definitions of differential calculus and connections on them. In Section \ref{section3}, we show that if the bimodule $ \E \tensora \E $ ( $ \E $ being the space of one forms coming from a differential calculus ) admits a splitting into symmetric and antisymmetric $2$-tensors, then $ \E $ admits a torsionless connection which is canonically related with the Grassmann connection. As a result, we have a symmetrization map $ \sigma $ on the bimodule $ \E \tensora \E $.    This allows us to define the notion of a pseudo-Riemannian metric and study its properties in Section \ref{section4}. Consequently, in Section \ref{section5}, we define the metric-compatibility condition of a connection on the center of the module $\E$ and prove a Koszul formula for a torsionless and metric-compatible (on the center) connection for bilinear pseudo-Riemannian metric. Under an additional assumption (see Theorem \ref{existenceuniqueness}) made in  Section \ref{section6}, we prove the existence and uniqueness of a torsionless and metric-compatible (on the whole of $\E$) connection as an application of the Koszul formula proved in Section \ref{section5}. In Section \ref{section7}, as a genuine application of the Koszul formula, we prove that our Levi-Civita connection is indeed a bimodule connection. Finally, in Section \ref{section8}, we construct a spectral triple ( see \cite[Def.~2, page 546]{connes} ) for the fuzzy sphere and prove the existence of the Levi-Civita connection and compute the scalar curvature for a canonical    bilinear pseudo-Riemannian metric.

We  fix some notations which we will follow. Throughout the article, $ \mathcal{A} $ will denote a complex algebra and $ \Acenter $ will denote its center. The tensor product over the complex numbers $ \IC $ is denoted by $ \tensorc $ while the notation $\tensora$ will denote the tensor product over the algebra $ \A$. For a subset $ S $ of a right $ \A $-module $ \E$, $S\A $ will denote 
its right $\A$-linear span: $ S\A = {\rm span} \{ sa: s \in S, ~ a \in \A  \}$.  
We will say that a subset $ S  $ of a right $ \A $-module $ \E $ is right $ \A $-total in $ \E $ if   the right $\A$-linear span of $ S  $ equals $ \mathcal{E}$.

For $\A-\A$-bimodules $ \E $ and $ \F, $ the symbol $ \Hom_\A  ( \E, \F ) $ will denote the set of all right  $ \A $-linear maps from $ \E $ to $ \F $. Similarly, $ {}_\A \Hom ( \E, \F ) $ will denote the set of all left $\A$-linear maps from $\E$ to $\F.$ In particular, we will  use the shorthand notation $ \E^* = \Hom_\A ( \E, \A )$.

For $\A-\A$-bimodules $ \mathcal{F} $ and $ \mathcal{F^\prime}, $ let us spell out the left and right $ \A $-module structures for   $ {\rm Hom}_\A ( \mathcal{F}, \mathcal{F^\prime} ) $ and $ {}_\A \Hom ( \F, \F^\prime ) $. 

The bimodule multiplications on $ \Hom_\A ( \F, \F^\prime ) $ and on $ {}_\A \Hom ( \F, \F^\prime ) $ are respectively:
\begin{align}
( a. T ) ( f ) &= a T ( f ) \in \F^\prime, \quad T. a ( f ) = T ( a f ),  \qquad a \in \A, f \in \F, T \in \Hom_\A ( \F, \F^\prime ); \nn \\
~ \nn \\
( a. T ) ( f ) &= T ( f a ), \quad ( T. a ) ( f ) = T ( f ). a, \qquad a \in \A, f \in \F, T \in {}_\A \Hom ( \F, \F^\prime ).  \label{22ndjuly20191}
\end{align}
 
\section{Differential calculus and connections on one-forms} \label{13thjuly20191}

As already mentioned, in the context of noncommutative geometry, it is more natural to work on the level of differential forms and that is what we do here. In this section, we recall the definition of connections on the space of one-forms coming from a differential calculus. This is followed by the notion of torsion. Then we define the notion of pseudo-Riemannian metric and compatibility of a connection on one forms under some assumptions on the differential calculus and the pseudo-Riemannian metric. 
			
			\bdfn Suppose $\A$ is an algebra over $\IC.$ A differential calculus on $\A$ is a pair $ (  \Omega ( \A ), d   ) $ such that the following conditions hold:
\begin{itemize}		 
		  \item[1.]  $\Omega ( \A )$ is graded: $ \Omega ( \A ) = \oplus_{ j \geq 0 } \Omega^j ( \A ), $ where $ \Omega^0 ( \A ) = \A $ and $ \Omega^j ( \A ) $ are $\A-\A$-bimodules. Thus, $ \Omega ( \A ) $ is an $ \A - \A $-bimodule.
			
			\item[2.] We have a bimodule map $ \wedge: \Omega ( \A ) \tensora \Omega ( \A ) \rightarrow \Omega ( \A ) $ such that 
			   $$ \wedge ( \Omega^j ( \A ) \tensora \Omega^k ( \A )  ) \subseteq \Omega^{j + k} ( \A )  .$$
				
		  \item[3.] We have a map $ d: \Omega^j ( \A ) \rightarrow \Omega^{j + 1} ( \A ) $ such that 
				$$ d^2 = 0 \qquad \textup{and} \qquad ( \omega \wedge \eta ) = d \omega \wedge \eta + ( - 1 )^{{\rm deg} ( \omega ) } \omega \wedge d \eta  .$$
				
			\item[4] $ \Omega^j ( \A ) $ is the right $ \A $-linear span of elements of the form $ d a_0 \wedge d a_1 \wedge \cdots \wedge d a_{j - 1} $.
	\end{itemize}
	  \edfn

\begin{assu}
\textup{		
Throughout the present paper, the notation $ \E $ will stand for the space of one-forms $\oneform$ of a differential calculus. It will also be assumed that $ \E $ is a finitely generated projective right $  \A  $-module.
}
\end{assu}		
		
		\bdfn \label{rLr}
Let $ ( \Omega ( \A ), d ) $ be a differential calculus on $\A.$  A (right) connection on $\E: = \oneform$
is a ${\mathbb C}$-linear map $\nabla :\E  \rightarrow \E \tensora \oneform$ satisfying the Leibniz rule 
$$ \nabla(\omega a)=\nabla(\omega)a + \omega \tensora da$$
for all $ \omega \in \E, a \in \A $.
\edfn

The assumption that $\E$ is finitely generated and projective is crucial in the present paper. In fact, by Corollary 8.2 of \cite{cuntz}, a (right) module admits a connection if and only if it is projective. Let us recall the construction of the Grassmann connection $ \nabla^{Gr} $. 

Since $ \E $ is finitely generated and projective as a right $ \A $-module, there exists a natural number $ n $ and an idempotent $ p \in M_n ( \A ) $ such that $ p ( \A^n ) = \E $. If $ \{ e_j: j = 1, \cdots n \} $ is a basis of the free right $ \A $-module $ \A^n, $ then the elements $ \{ \Phi_j: = p ( e_j ) : j = 1, \cdots n \} $ form a ``frame" ( in the sense of Rieffel, \cite{rieffelgh} ) of $ \E $. In particular, $ {\rm Span}_{\mathbb{C}} \{ \Phi_j \} $ is right $\A$-total in $ \E $.
Let $ \eta $ be an element in $ \E $. Then there exist elements $ \{ a_j: j = 1, \cdots n \} $ in $ \A $ such that 
$ \eta = \sum\nolimits_j \Phi_j a_j  $ and   the Grassmann connection $ \nabla^{Gr} $ is defined to be:
$$ \nabla^{Gr} ( \eta ) = \sum\nolimits_j \Phi_j \tensora d a_j  .$$
It is well-known that the set of all connections on $ \E $ is an affine space: any two right connections on $ \E $ differ by an element of $ \Hom_\A ( \E, \E \tensora \E ) $.  

\section{Existence of a torsionless connection} \label{section3}

We next recall the notion of the torsion of a connection and show the existence of a torsionless connection on a finitely generated projective module $ \E $ in the presence of a splitting of the right $ \A $-module $ \E \tensora \E $. 

\bdfn
	The torsion of a connection $ \nabla: \E \rightarrow \E \tensora \E $ is the right $\A$-linear map  
	$$ T_\nabla:= \wedge \circ \nabla + d: \E \rightarrow \twoform  .$$
	A connection $\nabla$  is called torsion-less if $ T_\nabla = 0.$
	\edfn
	
	The torsion of the Grassmann connection $ \nabla^{Gr} $ defined in Section \ref{13thjuly20191} is non-zero. Indeed, if $ \eta = \sum\nolimits_j \Phi_j a_j $ as in Section \ref{13thjuly20191}, then  
	\begin{align*}
T_{\nabla^{Gr}} ( \eta ) & = \wedge \nabla^{Gr} ( \eta ) + d ( \eta ) \\
	&= \sum\nolimits_j \Phi_j \wedge d a_j + \sum\nolimits_j d ( \Phi_j ) a_j - \sum\nolimits_j \Phi_j \wedge d a_j\\
	&= \sum\nolimits_j d ( \Phi_i ) a_j.
	\end{align*}
	
	\brmrk
	Let $ M $ be a manifold and $\A$ the algebra $ C^\infty ( M ) $. Let us consider the classical differential calculus $ ( \Omega ( \A ), d ) $ where $d$ is the de-Rham differential and $ \Omega ( \A ) $ the usual space of forms. Thus, $\E$ is the usual space of one-forms. It can be easily checked that the differential calculus $ ( \Omega ( \A ), d ) $ is tame in the sense of Definition 2.2 of \cite{article4}. Thus, we can apply Proposition 5.1 of \cite{article4} to conclude that our definition coincides with the usual definition of a torsionless connection in the classical case.
	\ermrk

	\bthm \label{torsionless}
	Suppose the short exact sequence of right $\A$-modules  
	$$ 0 \rightarrow {\rm Ker} ( \wedge ) \rightarrow  \E \tensora \E  \rightarrow {\rm Ran}(\wedge) = \twoform $$
	splits. Then there exists a torsionless connection $ \nabla_0 $ on $ \E $.
	\ethm
	\noindent {\bf Proof:} By our assumption, there is a right $\A$-submodule $ \F $ of $ \E \tensora \E $ and a  right $ \A $-module isomorphism $ Q: \F \rightarrow \twoform $ such that $ Q ( \beta ) = \wedge ( \beta ) $ for all $ \beta \in \F.$ 
 We define $ \nabla_0: \E \rightarrow \E \tensora \E $ by the formula:
  $$ \nabla_0 = \nabla^{Gr} - Q^{-1} ( T_{\nabla^{Gr}}  ).  $$
	Then $ \nabla_0 $ is a connection since for all $ \eta \in \E $ and for all $ a \in \A, $ we have
		\begin{align*}
		 \nabla_0 ( \eta a ) &= \nabla^{Gr} ( \eta a ) - Q^{-1} ( T_{\nabla^{Gr}} ( \eta a )  ) = \nabla^{Gr} ( \eta ) a + \eta \tensora da - Q^{-1} ( T_{\nabla^{Gr}} ( \eta )  ) a\\
		&=  ( \nabla^{Gr} ( \eta ) - Q^{-1} (  T_{\nabla^{Gr}} ( \eta )  )    ) a + \eta \tensora da = \nabla_0 ( \eta ) a + \eta \tensora da,
		\end{align*}
		where we have used that fact that $  Q $ is a right $ \A $-linear map.
		
Finally, $ \nabla_0 $ is torsionless as
\begin{align*}
 \wedge \nabla_0 ( \eta ) + d \eta &=   \wedge \nabla^{Gr} ( \eta ) - \wedge Q^{-1} (  T_{\nabla^{Gr}} ( \eta )  ) + d \eta\\
&= T_{\nabla^{Gr}} ( \eta ) - T_{\nabla^{Gr}} ( \eta ) = 0.
\end{align*}
This finishes the proof of the theorem. \qed

\bdfn \label{16thjuly201923}
Suppose $ \E $ satisfies the hypothesis of Theorem \ref{torsionless}. We will denote by the symbol $ \Psym $ the idempotent in $ \Hom_\A ( \E \tensora \E, \E \tensora \E ) $ with image $ {\rm Ker} ( \wedge ) $ and kernel $ \F $. Moreover, $ \sigma $ will be the map 
$$ \sigma = 2 \Psym - 1  .$$
\edfn
Let us also note that $ \sigma^2 = \id_{\E \tensora \E}$. 
Thus, we have
$$ \E \tensora \E = {\rm Ker} ( \wedge ) \oplus \F $$
where $ {\rm Ker} ( \wedge ) = {\rm Ran} ( \Psym ) $ and $ \F = {\rm Ran} ( 1 - \Psym ) $. Also, $ \F $ is isomorphic to $ \twoform $ as right $ \A $-modules via a right $ \A $-linear isomorphism $ Q: \F \rightarrow \twoform $. In fact, $ Q = \wedge  $ on $ \F $.

	We will need to define the action of $ \E^* \tensora \E^*   $ on the space of two forms. For that, let us recall that the map $ \wedge $ is an isomorphism from $ \F = {\rm Ran} ( 1 - \Psym ) $ onto $ \Omega^2 $.
	
	\bdfn \label{29thjandefn}
	Suppose $ \phi, \psi $ are elements of $ \E^* $ and let $ W $ be an element of $ \Omega^2 $. We define
	$$ ( \phi \tensora \psi ) W = 2 ( \phi \tensora \psi ) \beta  .$$
	where $ \beta $ is the unique element in $ \F = {\rm Ran} ( 1 - \Psym ) $ such that $ W = \wedge \beta $.
	\edfn
Here the factor $2$ is just a normalization factor in the product $\wedge$.	
	Let us note the following consequence of the definition:
	
	\blmma \label{29thjanlemma1}
	Suppose $ \phi, \psi, W, \beta $ be as in Definition \ref{29thjandefn} and $ \gamma $ an element of $ \E \tensora \E $ such that $ \wedge \gamma = W $. Then
	$$ ( \phi \tensora \psi ) W  = 2 ( \phi \tensora \psi ) ( 1 - \Psym ) ( \gamma ).  $$
	\elmma
	\noindent {\bf Proof:}  Since $ \wedge ( \gamma - \beta ) = 0$, then $\gamma - \beta ~ \in ~ {\rm Ker} ( \wedge ) = {\rm Ran} ( \Psym )   $ so that 
	$$ ( 1 - \Psym ) ( \gamma - \beta ) = 0  .$$
	Therefore,
	\begin{align*}
	 2 ( \phi \tensora \psi ) ( 1 - \Psym ) ( \gamma ) - 2 ( \phi \tensora \psi ) ( \beta ) &= 2 ( \phi \tensora \psi ) [ ( 1 - \Psym ) ( \gamma ) - ( 1 - \Psym ) ( \beta ) ]\\
	 &= 2 ( \phi \tensora \psi ) ( 1 - \Psym ) ( \gamma - \beta ) = 0.
	\end{align*}
	This proves the lemma.
		\qed

\section{Pseudo-Riemannian metrics on centered bimodules} \label{section4}
We now recall the notion of metric on a bimodule and work out some additional properties on a class of bimodules that we shall use in the rest of the paper. 
\bdfn \label{metricdefn}
Suppose $\E$ is an $\A - \A $-bimodule satisfying the hypothesis of Theorem \ref{torsionless} 
and let $\sigma$ the corresponding map as defined in Definition \ref{16thjuly201923}. 
  A pseudo-Riemannian metric $ g $ on $ \E $ is
 an element of $ {\rm Hom}_{\A} ( \E \tensora \E, \A ) $ such that
\begin{itemize} 
\item[(i)] $g$ is symmetric, that is $g \sigma = g$. 

\item[(ii)] The map $ \E \rightarrow \E^*, ~ e \mapsto g ( e \tensora - ) $
 is an isomorphism of right $ \A$-modules.
\end{itemize} 
We say that a pseudo-Riemannian metric $g$ is a pseudo-Riemannian bilinear metric if, in addition, $ g $ is an  $ \A - \A$ bimodule map.
In this case, the map $ g ( e \tensora - ) $ is bilinear as well.
\edfn

When $M$ is a manifold and $\A = C^\infty ( M )$ is the algebra of smooth complex valued functions on $M,$ a metric $g$ is a choice of a smooth positive definite symmetric bilinear form on the tangent ( or cotangent ) bundle. There are two equivalent ways to extend this map to the complexified spaces. The first way is to extend $g$ as a sesquilinear pairing on the module of one forms and thus linear in one variable and conjugate-linear in the other variable. In this paper, we have taken the second path, namely, we extend $g$ as a complex bilinear form, i.e, a  $C^\infty ( M ) $-bilinear map on $ \Omega^1 ( M ) \otimes_{C^\infty ( M )} \Omega^1 ( M )$.

In the noncommutative case, these two approaches do not remain equivalent. Since we do not deal with $\ast$-algebras in this paper we do not need any other compatibility with the $\ast$-structure as in \cite[Sect.~3]{pseudo}, 
or \cite[Def.~8.30]{beggsmajidbook}, or \cite{frolich}. We also mention the work \cite{heckenberger} which assumes the left $\A$-linearity of $g$ ( as opposed to right $\A$-linearity ) and derives the existence and uniqueness of Levi-Civita connections on three families of Hopf algebras. In \cite[ App.~B]{heckenberger} there is an example where the uniqueness of the Levi-Civita connection is lost when working with a sesquilinear metric.   Our main result ( Theorem \ref{existenceuniqueness} ) assumes both left and right $\A$-linearity of $g$, a condition which is restrictive. The conformally deformed metrics studied after the work of Connes, Moscovici, are only right $\A$-linear. The machinery developed in the present paper ( as well as \cite{article1} ) acts as a stepping stone for the existence and uniqueness ( for which we refer to \cite{article2} ) of Levi-Civita connections  for right $\A$-linear pseudo-Riemannian metrics on a differential calculus, satisfying the hypotheses of Theorem \ref{existenceuniqueness}.

In the present paper we shall be interested in a particular kind of bimodules that are called centered. 
Now, the center of  an $ \A - \A $-bimodule $ \mathcal{E} $ is defined to be the set 
$$ \mathcal{Z} ( \mathcal{E} ) = \{ e \in \mathcal{E}: e a = a e ~ \forall ~ a ~ \in \A  \}  .$$
 It is easy to see that $ \mathcal{Z} ( \E ) $ is a $ \mathcal{Z} ( \A ) $-bimodule.  The bimodule $ \mathcal{E} $ is called centered if $ \mathcal{Z} ( \mathcal{E} ) $ is right $ \A $-total in $ \E $, that is,  the right $\A$-linear span of $ \mathcal{Z} ( \mathcal{E} ) $ equals $ \mathcal{E}$.

Let us clarify that the property of being a centered bimodules is actually stronger than being a central bimodule in the sense of \cite{dubois} from which we have the following:
\bdfn
Suppose $ \A $ is a unital algebra and $\F$ is an $\A$-$\A$-bimodule. Then $\F$ is called a central bimodule if $ e. a = a. e  $ for all $e$ in $\E$ and for all $a$ in $ \Acenter $.
\edfn
It is easy to see that a centered bimodule is a central bimodule. Indeed, if $ \E $ is a centered bimodule, then for any $e \in \E$ there exists a natural number $ n, $ elements $ f_j \in \Ecenter $ and $ b_j \in \A $ such that 
$ e = \sum\nolimits_j f_j b_j $.
Then 
\begin{equation} \label{16thjuly201921}a. e = \sum\nolimits_j a . f_j. b_j = \sum\nolimits_j f_j. a. b_j = \sum\nolimits_j f_j. b_j. a. = e. a  
\end{equation}
for all $e \in \E$ and $a \in \Acenter$. Thus, $ \E $ is central.

\blmma \label{lemma2}
Suppose  $g$ is a pseudo-Riemannian metric on a centered $\A - \A $-bimodule  $\E$ satisfying the hypothesis of Theorem \ref{torsionless} with $\sigma$ the corresponding map as defined in Definition \ref{16thjuly201923}. 
 Moreover, assume
 $$ \sigma ( \omega \tensora \eta ) = \eta \tensora \omega \qquad \forall ~ \omega, \eta \in \Ecenter  .$$
Then we have the following: 
\begin{itemize}
\item[1.] If either  $ \omega $ or $ \eta $ belongs to $ \Ecenter, $ then $ \sigma ( \omega \tensora \eta ) = \eta \tensora \omega $.
\item[2.]If   either of $ \omega $  or $ \eta $ belongs to $ \Ecenter,$ then    
\begin{equation} \label{gsigmaisg}
 g ( \omega \tensora \eta ) = g ( \eta \tensora \omega ). 
\end{equation}
\item[3.] If $ g $ is a  pseudo-Riemannian  bilinear metric, then $ g ( \omega \tensora \eta ) \in \mathcal{Z} ( \A ) $ if both $ \omega $ and $ \eta  $ are in $ \mathcal{Z} ( \E ).$
\item[4.] If $ a $ is an element of $ \Acenter, $ then $ da \in \Ecenter $. In particular,
if $ \omega, \eta \in \Ecenter $ and $ g $ is a pseudo-Riemannian bilinear metric, then 
\begin{equation} \label{22ndjuly20192} d \, g ( \omega \tensora \eta ) \in \Ecenter. \end{equation}
\end{itemize}
\elmma
\noindent {\bf Proof:} The first three assertions were proved in Lemma 2.8 of \cite{article1}. As for 4., for all $ b \in \A, $ we have
	$$ 0 = d ( a. b ) - d ( b. a ) = ( da. b + a. db    ) - ( db. a + b. da   ) = da. b - b. da $$
	and we have used \eqref{16thjuly201921}. This proves that $ da \in \Ecenter $.
	
	By part 3., $ g ( \omega \tensora \eta ) \in \Acenter $ and hence $ d ( g ( \omega \tensora \eta )  ) \in \Ecenter $. This proves \eqref{22ndjuly20192}. 
	\qed

\brmrk \label{21stjuly201910}
It is easy to see that for $ \xi \in \E, $ the condition $ g ( \theta \tensora \xi ) = 0 $ for all $ \theta \in \Ecenter $ implies that $ \xi = 0 $.
\ermrk

\section{The Koszul formula on one-forms on the center of the module} \label{section5}

Throughout this section, we will assume that $ g $ is a bilinear pseudo-Riemannian metric so that the assertions of Lemma \ref{lemma2} are valid. Moreover, $ \nabla_0 $ will denote the torsionless connection of Theorem \ref{torsionless}.

\subsection{Metric compatibility of a connection on the center}

Let $ g $ be a pseudo-Riemannian metric on $\E$ and $ \nabla $ a connection on $ \E $. It can be checked that the map
\begin{align*}
& \Ecenter \otimes_{\Acenter} \Ecenter \rightarrow \E,  \\ 
& \omega \otimes_{\Acenter} \eta \, \mapsto \, 
(g \tensora \id ) \Big\{ \big[\sigma_{23}(\nabla(\omega)\tensora \eta )\big] + \omega \tensora \nabla ( \eta ) \Big\}
\end{align*}
 is well defined. 
 Indeed, for $ a \in \Acenter, $ 
we have $\sigma_{23}(\omega\tensora da\tensora \eta)=\omega\tensora \eta\tensora da$, $ a \eta = \eta a $ and $ a \nabla ( \eta ) = \nabla ( \eta ) a $ ( by \eqref{16thjuly201921} ). 
Using these, we get
\begin{align*}
\sigma_{23} (\nabla(\omega)a\tensora \eta & + \omega \tensora da \tensora \eta ) +  \omega a \tensora \nabla ( \eta )  \\
 & = \sigma_{23}(\nabla(\omega) \tensora a \eta + \omega \tensora da \tensora \eta ) + \omega \tensora a \nabla ( \eta )  \\ 
& = \sigma_{23}(\nabla(\omega)\tensora a \eta ) +  \sigma_{23} ( \omega \tensora d a \tensora \eta    ) + \omega \tensora \nabla( \eta ) a \\
& = \sigma_{23}(\nabla(\omega) \tensora a \eta  ) + \omega \tensora \eta \tensora da + \omega \tensora \nabla ( \eta a ) - \omega \tensora \eta \tensora da \\
& = \sigma_{23}(\nabla(\omega)\tensora a \eta )\big] + \omega \tensora \nabla ( a \eta ) 
\end{align*}
This proves the well-definedness. 
 
\bdfn \label{compatibilitycenter}
We say that  a connection $\nabla$ on $\E$ is compatible with $g$ on $ \Ecenter $ if for all $ \omega, \eta \in \Ecenter, $ the following equation holds:
$$
(g \tensora {\rm id} ) \Big\{ \big[\sigma_{23}(\nabla(\omega)\tensora \eta )\big] + \omega \tensora \nabla ( \eta ) \Big\}
 = d ( g ( \omega \otimes_{\Acenter} \eta ) )  .
$$
\edfn

\brmrk  Given that $\nabla(\omega) \tensora \eta \in \E \tensora \Omega^1 \tensora \E$ and 
$\omega \tensora \nabla ( \eta ) \in \E \tensora \E \tensora \Omega^1, $ we can write the metric compatibility condition on the center in another equivalent way. 
 The condition is:
$$( {\rm id} \tensora g )\sigma_{12}(\nabla(\omega) \tensora \eta ) + ( g \tensora {\rm id}  ) ( \omega \tensora \nabla ( \eta ) ) 
= d ( g ( \omega \tensora \eta ) ) ~ \forall \omega, \eta \in \Ecenter  .$$
Here, we have used the fact that $ g $ is bilinear so that $ ( {\rm id} \tensora g  ) $ is well-defined.   The proof of this equation is a straightforward application of the facts that $ \sigma ( e \tensora f ) = f \tensora e $ if either $e$ or $ f $ belong to $\Ecenter$ and that we can write  $\nabla ( e ) = \sum\nolimits_{j} f_j \tensora \omega_j $ with $ \omega_j \in \Ecenter $
(see Part 3. in Lemma \ref{20thaugust20191}).
      \ermrk

\subsection{Consequences of the zero-torsion condition}

\begin{nota}\label{notasweed}
\textup{
For $ \nabla $ a torsionless connection on a centered bimodule $\E$ we shall use Sweedler-like  notation to write 
\begin{equation}\label{sweedler}
\nabla ( \omega ) = \sum\nolimits_j \omegazero{}_j \tensora \omegaone{}_j =: \omegazero \tensora \omegaone 
\end{equation}
Likewise, for the torsionless connection $ \nabla_0 $ of Theorem \ref{torsionless}, we will write
\begin{equation}\label{sweedler0}
\nabla_0 ( \omega ) = \sum\nolimits_j \zeroomega{}_j \tensora \oneomega{}_j =: \zeroomega \tensora \oneomega.  
\end{equation}
One can always take that both $\omegaone $  and $ \oneomega $ belongs to $\Ecenter$
( see Parts 3. and 4. of Lemma \ref{20thaugust20191} ). 
}
\end{nota}

	\blmma \label{13thjuly20192}
	Suppose	$ \nabla $ is a torsionless connection on a centered bimodule $\E$ and $\nabla_0 $
	the torsionless connection $ \nabla_0 $ of Theorem \ref{torsionless} for which we use the notation \ref{notasweed} with both $\omegaone $  and $ \oneomega $ in $\Ecenter$. 
	
Then for all $ \omega, \eta, \theta $ in $ \E, $ the following equations hold:
\begin{multline*} g ( \eta \tensora \omegazero ) \, g ( \theta \tensora \omegaone ) - g ( \eta \tensora \omegaone ) \, g ( \theta \tensora \omegazero )  \\ =  g ( \eta \tensora \zeroomega )  \, g ( \theta \tensora  \oneomega ) - g ( \eta \tensora \oneomega )  \, g ( \theta \tensora  \zeroomega ), 
\end{multline*}
\begin{multline*}
	 g ( \omega \tensora \etazero ) g ( \theta \tensora \etaone ) + g ( \omega \tensora \etaone  ) g ( \theta \tensora \etazero  )  \\ = 2 g ( \omega \tensora \etazero ) g ( \theta \tensora \etaone )
	-  g ( \omega \tensora \zeroeta  ) g ( \theta \tensora \oneeta )\\
	+ g ( \omega \tensora  \oneeta ) g ( \theta  \tensora \zeroeta  ),
\end{multline*}
\begin{multline*} g ( \eta \tensora \thetazero ) g  ( \omega \tensora \thetaone ) - g ( \eta \tensora \thetaone ) g ( \omega \tensora \thetazero ) \\ =  g ( \eta \tensora \zerotheta ) g ( \omega \tensora \onetheta ) - g ( \eta \tensora \onetheta ) g ( \omega \tensora \zerotheta ).  
\end{multline*}
	\elmma
		\noindent {\bf Proof:} We will use Lemma \ref{29thjanlemma1}. Suppose $ \omega \in \E $. Since $ \nabla $ and $ \nabla_0 $ are both torsionless, we have
			$$ 0 = \wedge ( \nabla ( \omega ) - \nabla_0 ( \omega ) )  .$$
	Using Lemma \ref{29thjanlemma1},  this implies that for all $ \phi, \psi \in \E^*, $ we have
									\begin{align*}
 0 &= ( \phi \tensora \psi  ) \, \wedge \big( \nabla ( \omega )   -  \nabla_0 ( \omega ) \big)\\
 &= 2 ( \phi \tensora \psi  ) \, ( 1 - \Psym ) \, \big( \nabla ( \omega )   -  \nabla_0 ( \omega ) \big)\\
  &= ( \phi \tensora \psi  ) \, ( 1 - \sigma ) \, \big( \nabla ( \omega )   -  \nabla_0 ( \omega ) \big)\\
					&=  ( \phi \tensora \psi ) ( \omegazero \tensora \omegaone - \omegaone \tensora \omegazero ) - ( \phi \tensora \psi  ) ( \zeroomega \tensora \oneomega - \oneomega \tensora \zeroomega ).  
				\end{align*}
				
				Therefore, we obtain
				$$ ( \phi \tensora \psi ) ( \omegazero \tensora \omegaone - \omegaone \tensora \omegazero ) = ( \phi \tensora \psi  ) ( \zeroomega \tensora \oneomega - \oneomega \tensora \zeroomega ).$$
				Putting $ \phi = g ( \eta \tensora - )$ and $\psi = g ( \theta \tensora - )$, we obtain the first equation. The other two equations are obtained similarly.
				\qed
							
	\subsection{A Koszul formula for the Levi-Civita connection on the center}
							
	Let $ \nabla $ be a torsionless connection on $ \E $ which is compatible with a pseudo-Riemannian metric on $ \Ecenter $ as in Definition \ref{compatibilitycenter}. We will still use the Sweedler type notations $ \nabla ( \omega ) = \omegazero \tensora \omegaone$ and $\nabla_0 ( \omega ) = \zeroomega \tensora \oneomega, $
		where $ \omegaone$ and $\oneomega $ belong to $\Ecenter.$

\bthm \label{koszul21stjuly}
Let $ \omega, \eta, \theta \in \Ecenter  $ and $ \nabla $ is a torsionless connection on $ \E $ which is compatible with a bilinear pseudo-Riemannian metric $ g $. Then the following formula holds:
\begin{align} \label{17thjuly20191}
  2  g ( \omega \tensora \etazero ) & g ( \theta \tensora \etaone ) \nn \\ & =  g ( \omega \tensora dg ( \eta \tensora \theta ) )  - g ( \eta \tensora dg ( \theta \tensora \omega ) ) + g ( \theta \tensora dg ( \omega \tensora \eta ) ) \nonumber \\
  & \quad - g ( \eta \tensora \zeroomega ) \,\, g ( \theta \tensora \oneomega ) + g ( \eta \tensora \oneomega ) \,\, g ( \theta \tensora \zeroomega ) \nonumber \\
	& \quad + g ( \omega \tensora \zeroeta ) \,\, g ( \theta \tensora \oneeta ) - g ( \omega \tensora \oneeta ) \,\, g ( \theta \tensora \zeroeta ) \nonumber \\
	& \quad - g ( \eta \tensora \zerotheta ) \,\, g ( \omega \tensora \onetheta ) + g ( \eta \tensora \onetheta ) \,\, g ( \omega \tensora \zerotheta ).
\end{align}
\ethm
\noindent {\bf Proof:} Since $ \nabla $ is compatible with $ g $ on $ \Ecenter, $ ( Definition \ref{compatibilitycenter}
 ) we have
$$ g ( \omegazero \tensora \eta  ) \omegaone + g ( \omega \tensora \etazero ) \etaone = d ( g ( \omega \tensora \eta ) )  .$$
In turn, since $ \omegaone, \etaone \in \Ecenter, $ this implies
$$ \omegaone g ( \omegazero \tensora \eta  ) + \etaone g ( \omega \tensora \etazero ) = d ( g ( \omega \tensora \eta ) ).  $$
Applying $ g ( \theta \tensora - ) $ to the above equation, we get
\begin{equation} \label{13thjuly20194} g ( \theta \tensora \omegaone ) g (  \omegazero \tensora \eta ) + g ( \theta \tensora \etaone ) g ( \omega \tensora \etazero ) = g ( \theta \tensora dg ( \omega \tensora \eta ) ). \end{equation}
Replacing $ \omega, \eta, \theta $ by $ \eta, \theta, \omega  $ respectively in \eqref{13thjuly20194}, we get
\begin{equation} \label{13thjuly20195} g ( \omega \tensora \etaone ) g (  \etazero \tensora \theta ) + g ( \omega \tensora \thetaone ) g ( \eta \tensora \thetazero ) = g ( \omega \tensora dg ( \eta \tensora \theta ) ). \end{equation}
Replacing $ \omega, \eta, \theta $ by $ \theta, \omega, \eta  $ respectively in \eqref{13thjuly20194}, we get
\begin{equation} \label{13thjuly20196} g ( \eta \tensora \thetaone ) g (  \thetazero \tensora \omega ) + g ( \eta \tensora \omegaone ) g ( \theta \tensora \omegazero ) = g ( \eta \tensora dg ( \theta \tensora \omega ) ). \end{equation}
By \eqref{13thjuly20194} + \eqref{13thjuly20195} - \eqref{13thjuly20196}, 
using $ g \sigma = g$ and the fact that $ g (\theta \tensora \omegaone )$, $g ( \theta \tensora \etaone )$ and
$g ( \omega \tensora \thetaone )$, $g ( \eta \tensora \thetaone )$ all belong to $\Acenter $ ( Lemma \ref{lemma2} ),
we obtain
\begin{align*}
&  g ( \eta \tensora \omegazero  ) g ( \theta \tensora \omegaone  ) - g ( \eta \tensora \omegaone ) g ( \theta \tensora \omegazero  ) ] 
+ [  g ( \omega \tensora \etazero  ) g ( \theta \tensora \etaone  )\\
 & +  g ( \omega \tensora \etaone   ) g ( \theta \tensora \etazero  ) ] + [ g  ( \eta \tensora \thetazero  ) g ( \omega \tensora \thetaone  ) - g ( \eta \tensora \thetaone ) g ( \omega \tensora \thetazero  )   ]\\
 & \quad = g ( \theta \tensora dg ( \omega \tensora \eta ) ) +  g ( \omega \tensora dg ( \eta \tensora \theta ) ) - g ( \eta \tensora dg ( \theta \tensora \omega ) ).
\end{align*}
By Lemma \ref{13thjuly20192}, the left hand side of the above equation coincides with
\begin{align*}
  2  g ( \omega \tensora \etazero ) g ( \theta \tensora \etaone ) & + g ( \eta \tensora \zeroomega ) g ( \theta \tensora \oneomega ) + g ( \eta \tensora \oneomega ) g ( \theta \tensora \zeroomega )\\
	& -  g ( \omega \tensora \zeroeta ) g ( \theta \tensora \oneeta ) + g ( \omega \tensora \oneeta ) g ( \theta \tensora \zeroeta )\\
	& +  g ( \eta \tensora \zerotheta ) g ( \omega \tensora \onetheta ) + g ( \eta \tensora \onetheta ) g ( \omega \tensora \zerotheta ).
\end{align*}
This proves the theorem.
\qed

Let us end this section by rewriting the Koszul formula in an alternative form.
With $ \nabla_0 $ the torsionless connection of Theorem \ref{torsionless}, for all $ \omega, \eta, \theta \in \Ecenter, $ 
the following identity holds:
\begin{align} \label{16thjuly2019koszul}
2 ( g ( \omega \tensora - ) & \tensora g ( \theta \tensora - )   ) ( \nabla ( \eta ) ) 
\nn \\ &= g ( \omega \tensora dg ( \eta \tensora \theta ) ) - g ( \eta \tensora dg ( \theta \tensora \omega ) ) + g ( \theta \tensora dg ( \omega \tensora \eta ) )  \nonumber \\
& \quad - ( g ( \eta \tensora - ) \tensora g ( \theta \tensora - ) ) ( 1 - \sigma ) \nabla_0 ( \omega ) \nonumber 
\nn \\ 
 & \quad  +   (  g ( \omega \tensora - ) \tensora g ( \theta \tensora - ) ) ( 1 - \sigma ) \nabla_0 ( \eta ) \nn \\
& \quad - ( g ( \eta \tensora - ) \tensora g ( \omega \tensora - )   ) ( 1 - \sigma ) \nabla_0 ( \theta ).  
\end{align}

\section{Levi-Civita connection on a class of centered bimodules: existence} \label{section6}

The goal of this section is to prove the existence and uniqueness of Levi-Civita connections for a class of centered bimodules. Let us now state the main result precisely.

\bthm \label{existenceuniqueness}
Suppose $ ( \Omega^\cdot ( \A ), d ) $ is a differential calculus on  $ \A $ such that the following conditions are satisfied:
\begin{itemize}
\item[1.] The space of one-forms $ \oneform:=\E = \Ecenter \otimes_{\Acenter} \A $.

\item[2.]  We have $ \E \tensora \E = {\rm Ker} ( \wedge ) \oplus \F  $ where $ Q: \F \rightarrow \twoform  $ is  a right $\A$-linear isomorphism  as in Theorem \ref{torsionless}. 

\item[3.] The map $ \sigma $ in Definition \ref{16thjuly201923} satisfies $ \sigma ( \omega \tensora \eta ) = \eta \tensora \omega $ for all $ \omega, \eta \in \Ecenter $.
\end{itemize}

If $ g $ is a pseudo-Riemannian bilinear metric on $ \E, $ then there exists a unique connection on $ \E $ which is torsionless and compatible with $ g $ on $ \Ecenter $ as in Definition \ref{compatibilitycenter}.
\ethm

Throughout this section, we will work under the assumptions of Theorem \ref{existenceuniqueness}.

\subsection{A class of centered bimodules and examples} 

For proving the existence and uniqueness of the Levi-Civita connection, we will need some additional results which we spell out here. 
  
\blmma \label{20thaugust20191}
Let $\E = \Ecenter \otimes_{\Acenter} \A $. Then

1 $. \E $ is centered.

2 $. \Z(\E) $ is both left and right $ \A $-total in $ \E $.

3.  We have a right identification,
 $$
 \E \otimes_\A  \E \simeq \E \otimes_{\Z(\A)} \Z(\E)  .
 $$

4. as well as a left one,
 $$ 
 \E \otimes_\A  \E \simeq \Z(\E)  \otimes_{\Z(\A)} \E  .
 $$
 
5. The collection $ \{ \omega \tensora \eta : \omega, \eta \in \Ecenter \} $ is right $ \A $-total in $ \E \tensora \E $ .

6. The collection $ \{ \omega \tensora \eta : \omega, \eta \in \Ecenter \} $ is left $\A$-total in $\E \tensora \E $

\elmma
\noindent {\bf Proof:} Most of the statements in the lemma follow from the Proposition 2.4 of \cite{article1}. Indeed, that  proposition 
implies that the equality $ \E = \Ecenter \otimes_{\Acenter} \A $ forces $ \E $ to be centered. Moreover, we have right $\A$-module isomorphisms 
\begin{equation} \label{20thagust20192} \E \cong \A \otimes_{\Acenter} \Ecenter \cong \Ecenter \otimes_{\Acenter} \A \end{equation}
 via the canonical multiplication maps. The isomorphisms in particular imply that $ \Ecenter $ is both left and right $\A$-total in $ \E $. Thus, we have proved the parts 1,2, 3 and 4 of the lemma.

Thus we are left to prove the last two assertions. We only prove the last but one since the proof of the last one is similar.
If $ e,f $ are elements of $ \E, $ then \eqref{20thagust20192} allows us to conclude that there exist elements $ \omega_j, \eta_k, \eta_l \in \Ecenter $ and $  a_j, b_k, c_l \in \A$ such that
$$ e = \sum\nolimits_j a_j \omega_j, \quad f = \sum\nolimits_k b_k \eta_k, \quad \omega_j b_k = \sum\nolimits_l c_l \theta_l $$
Therefore, 
$$ e \tensora f = \sum\nolimits_{j, k} a_j \omega_j \tensora b_k \eta_k = \sum\nolimits_{j,k} a_j ( \omega_j b_k ) \tensora \eta_k = \sum\nolimits_{j,k,l} a_j c_l \theta_l \tensora \eta_k $$
which belongs to  the left $ \A $-linear span of $ \{ \omega \tensora \eta: \omega, \eta \in \Ecenter  \} $. This proves the fourth assertion.
\qed

 We will also need the following results about the maps $ \sigma, \Psym $ and $ Q $.
\bppsn
Suppose $ \E $ is an $ \A - \A $-bimodule satisfying the hypothesis of Theorem \ref{existenceuniqueness}. Then we have the following:

1. The map $ \sigma $ and hence $ \Psym $ is $\A - \A $-bilinear.

2. The map $ Q: \F \rightarrow \twoform $ is an $\A - \A $-bilinear isomorphism.

\eppsn
\noindent {\bf Proof:} Since  $ \E $ centered by part 1. of Lemma \ref{20thaugust20191}, Theorem 6.10 of \cite{Skd} implies there exists a unique $\A-\A$-bimodule isomorphism $\sigma^{{\rm can}}: \E \tensora \E \rightarrow \E \tensora \E$ such that $\sigma^{{\rm can}}(\omega \tensora \eta)=\eta \tensora \omega$ for all $\omega,\eta \in \Ecenter$. Thus, it suffices to show that $ \sigma = \sigma^{{\rm can}}.$ This follows once we show that $ \sigma $ is left $ \A $-linear. Let $ \omega, \eta \in \Ecenter $ and $ a, b \in \A $.
Then 
$$ a \sigma ( \omega \tensora \eta b ) = a \sigma ( \omega \tensora \eta ) b = a \eta \tensora \omega b = \eta \tensora \omega a b = \sigma ( a \omega \tensora \eta b )  .$$ 
Since $ \{ \omega \tensora \eta : \omega, \eta \in \Ecenter \} $ is right $ \A $-total in $ \E \tensora \E $ 
( part 5. of Lemma \ref{20thaugust20191} ), this proves that $ \sigma $ is left $ \A $-linear and completes 
the proof of part 1.

In particular, this implies that $ \F = {\rm Ran} ( 1 - \Psym ) $ is an $ \A - \A $-bimodule. Since $ \wedge = Q $ on $ \F $ and $ \wedge $ is bilinear, so is  $ Q $, thus establishing part 2.
\qed

\blmma \label{16thjuly20192}
We have the following:

a. Let $ \omega \in \Ecenter $. Then $ d \omega \in \mathcal{Z} ( \twoform ) $.

b. The element $ ( 1 - \sigma ) \nabla_0 ( \omega ) \in \mathcal{Z} ( \E \tensora \E ) $ if $ \omega \in \Ecenter $. 
\elmma
\noindent {\bf Proof:} For part a., we start by observing that 
\begin{align*}
\omega \wedge da + da \wedge \omega & = \wedge ( \omega \tensora da + da \tensora \omega )\\
& = \wedge ( 1 - \Psym ) ( \omega \tensora da + da \tensora \omega   ) = \tfrac{1}{2} ( 1 - \sigma ) ( \omega \tensora da + da \tensora \omega )\\ 
& = \tfrac{1}{2} \wedge ( \omega \tensora da + da \tensora \omega - da \tensora \omega - \omega \tensora da  ) = 0.
\end{align*}
Hence, 
$$ 0 = d ( \omega . a ) - d ( a. \omega )  = (  d \omega . a - a. d \omega   ) - (  \omega \wedge da + da \wedge \omega   ) = d \omega . a - a. d \omega.
$$
This proves part a.

For part b: being $ \nabla_0 $ a torsionless connection, we have $ \wedge \nabla_0 ( \omega ) = - d \omega $. By applying the map $ ( 1 - \sigma ) Q^{-1}, $ we get
$$ ( 1 - \sigma ) \nabla_0 ( \omega ) = - ( 1 - \sigma ) Q^{-1} ( d \omega )  .$$
Since $ Q^{-1} ( d \omega ) \in {\rm Ran} ( 1 - \Psym ), $ we have
$$ - ( 1 - \sigma ) Q^{-1} ( d \omega ) = - 2 ( 1 - \Psym ) ( Q^{-1} d \omega  ) = - 2 Q^{-1} ( d \omega ) .$$
Now, by part a $.d \omega \in \mathcal{Z} ( \twoform )$ since $\omega \in \Ecenter $; then since 
the map $ Q $ is bilinear, for any $ a \in \A $, 
$$ 
a. Q^{-1} ( \omega ) = Q^{-1} ( a. \omega ) = Q^{-1} ( \omega. a ) = Q^{-1} ( \omega ).a 
$$
and thus 
$Q^{-1} ( d \omega )$ belongs to $ \mathcal{Z} ( \E \tensora \E ) $. 

\qed

\brmrk \label{22ndjuly20193}
More generally, if $ \F $ and $ \F^\prime $ are centered $\A-\A$-bimodules and $ T: \F \rightarrow \F^\prime $ be an $ \A - \A $-bilinear map. Then $ T ( \Ecenter ) \subseteq \mathcal{Z} ( \F ) $. Indeed, 
If $ \omega \in \Ecenter $ and $ a \in \A $ a simple computation yields $ a. T ( \omega ) = T ( a. \omega ) = T ( \omega. a ) = T ( \omega ).a $ and thus, $ T ( \omega ) \in \Ecenter $. 
\ermrk

\subsection{The existence and uniqueness of Levi-Civita connection}

\bppsn \label{vg2nondegenerate} 
 Let $ g $ be a pseudo-Riemannian bilinear metric on $\E$ and define a map
\begin{align*}
& g^{(2)} :(\E\tensora \E)\tensora (\E \tensora \E) \rightarrow \A, \\
& g^{(2)}((e\tensora f)\tensora (e'\tensora f')) = g(e \tensora g(f \tensora e') f').
 \end{align*}
 Then we have the following:
\begin{itemize}
\item[1.] the map $ : \E \tensora \E \rightarrow ( \E \tensora \E )^* $ defined by
$$  e\tensora f \mapsto g^{(2)}((e \tensora f) \tensora - )$$
 is $\A-\A$-bilinear and an isomorphism of right $ \A $ modules.

\item[2.] For an element $ \xi \in \E \tensora \E, $ the map $ g^{(2)} ( - \tensora \xi ) $ is an element of $ {}_\A \Hom ( \E \tensora \E, \A ) $. Moreover, the map from $ \E \tensora \E $ to $ {}_\A \Hom ( \E \tensora \E, \A ) $ defined by
$$ \xi \mapsto g^{(2)} ( - \tensora \xi ) $$
is an isomorphism of right $ \A $-modules.
\end{itemize}
\eppsn
\noindent {\bf Proof:} We only need to prove 2. since the assertion 1. was already proved in Proposition 3.7 of \cite{article1}. However, the assertion 2. follows exactly as in the proof of Proposition 3.7 of \cite{article1} using the bimodule structure of $ {}_\A \Hom ( \E \tensora \E, \A ) $ as spelled out in \ref{22ndjuly20191}. 
\qed
\brmrk \label{11thaugust2019}
From the proof of Proposition 3.7 of \cite{article1}, 
it is easy to see that for $ \omega, \eta \in \Ecenter $ and  $ e, f \in \E, $
$$ ( g ( \omega \tensora - ) \tensora g ( \eta \tensora - ) ) ( e \tensora f ) = g^{(2)} ( ( \eta \tensora \omega  ) \tensora ( e \tensora f )   )  .$$
\ermrk

\blmma \label{21stjuly20192}
For $ \omega, \eta, \theta \in \Ecenter, $ let us define $ \psi_{\omega, \theta} ( \eta ) $ by  the expression of the   right hand side of \eqref{17thjuly20191} in Theorem \ref{koszul21stjuly}. Thus we have, 
\begin{align}\label{21stjuly20197}
\psi_{\omega, \theta} ( \eta ) 
&=   g ( \omega \tensora dg ( \eta \tensora \theta ) ) - g ( \eta \tensora dg ( \theta \tensora \omega ) ) + g ( \theta \tensora dg ( \omega \tensora \eta ) )  \nn \\
& \quad - ( g ( \eta \tensora - ) \tensora g ( \theta \tensora - ) ) ( 1 - \sigma ) \nabla_0 ( \omega ) \nn \\ 
& \quad +   (  g ( \omega \tensora - ) \tensora g ( \theta \tensora - ) ) ( 1 - \sigma ) \nabla_0 ( \eta ) \nn \\
 & \quad - ( g ( \eta \tensora - ) \tensora g ( \omega \tensora - )   ) ( 1 - \sigma ) \nabla_0 ( \theta ).
\end{align}
Then, for $ \omega, \eta, \theta \in \Ecenter, $ the element  $ \psi_{\omega, \theta} ( \eta )  $ belongs to $ \Acenter $.
\elmma
\noindent {\bf Proof:} Let us observe that since $ \omega, \eta, \theta \in \Ecenter, $ and $ d g ( \eta \tensora \theta ), d g ( \theta \tensora \omega ), d g ( \omega \tensora \eta )   $ are in  $ \Acenter $, hence by Lemma \ref{lemma2}  
$ g ( \omega \tensora dg ( \eta \tensora \theta ) ), g ( \eta \tensora dg ( \theta \tensora \omega ) )   $ and $ g ( \theta \tensora dg ( \omega \tensora \eta ) ) $ are also in $ \Acenter $. 

\noindent
Moreover, by part b. of Lemma \ref{16thjuly20192}, the element $ ( 1 - \sigma ) \nabla_0 ( \omega ) \in \mathcal{Z} ( \E \tensora \E ) $.  Therefore, for all $ a \in \A, $
\begin{align*}
( g ( \eta \tensora - ) \tensora g ( \theta \tensora - ) ) ( 1 - \sigma ) \nabla_0 ( \omega ). a 
&= g^{(2)} ( \theta \tensora \eta \tensora -  ) ( a. ( 1 - \sigma ) \nabla_0 ( \omega ) )\\
&= ( g^{(2)} ( \theta \tensora \eta \tensora -  ) . a  ) (  ( 1 - \sigma ) \nabla_0 ( \omega ) )\\
&= a. g^{(2)} ( \theta \tensora \eta \tensora -  ) ( ( 1 - \sigma ) \nabla_0 ( \omega ) )\\
\end{align*} 
as the map $ g^{(2)} ( \theta \tensora \eta \tensora -  ) $ is bilinear by Proposition \ref{vg2nondegenerate}. Hence, the element 
$$ ( g ( \eta \tensora - ) \tensora g ( \theta \tensora - ) ) ( 1 - \sigma ) \nabla_0 ( \omega ) $$ 
belongs to $ \Acenter $. Similarly, it belong to $ \Acenter $ also the elements  
$$ ( g ( \omega \tensora - ) \tensora g ( \theta \tensora - ) ) ( 1 - \sigma ) \nabla_0 ( \eta ) \qquad \mbox{and} \qquad 
( g ( \eta \tensora - ) \tensora g ( \omega \tensora - ) ) ( 1 - \sigma ) \nabla_0 ( \theta )   .$$ 
This completes the proof that $ \psi_{\omega, \theta} ( \eta )  $ belongs to $ \Acenter $.
\qed

\blmma \label{20thjuly20191}
Let $ \eta \in \Ecenter $. There exists a right $ \Acenter $-linear map
 $$ \phi_\eta: \Ecenter \otimes_{\Acenter} \Ecenter \rightarrow \Acenter \qquad \textup{such that} \qquad \phi_\eta ( \theta  \otimes_{\Acenter}  \omega ) = \psi_{\omega, \theta} ( \eta ) $$
\elmma
\noindent {\bf Proof:} We need to check that $ \phi_\eta $ is well-defined, i.e, that for all $ a^\prime \in \Acenter, $
$$ \phi_\eta ( \theta a^\prime  \tensorc \omega ) = \phi_\eta ( \theta \tensorc a^\prime \omega ).$$
Since $ ( 1 - \sigma )  ( \theta \tensora d a^\prime ) = 0,  $ we get 
\begin{align*}
\phi_n ( \theta a^\prime \tensora \omega ) &= \psi_{\omega, \theta a^\prime} ( \eta )\\
 &= ( g ( \eta \tensora - ) \tensora g ( \omega \tensora - )  ) ( 1 - \sigma ) \nabla_0 ( \theta a^\prime )\\
 &= ( g ( \eta \tensora - ) \tensora g ( \omega \tensora - )  ) ( 1 - \sigma ) \nabla_0 ( \theta ) a^\prime \\
 &= \psi_{\omega,\theta} ( \eta ) a^\prime\\
&= a^\prime  \psi_{\omega,\theta} ( \eta ) \qquad {\rm ( } {\rm Lemma} ~ \ref{21stjuly20192}  {\rm )}\\ 
&= a^\prime ( g ( \eta \tensora - ) \tensora g ( \omega \tensora - )  ) ( 1 - \sigma ) \nabla_0 ( \theta ).
\end{align*}
Now, since $ \eta \in \Ecenter $ and $ a^\prime \in \Acenter, $
$$ a^\prime. g ( \eta \tensora - ) = g ( \eta \tensora - ). a^\prime $$ 
and so
\begin{align*}
\phi_\eta ( \theta a^\prime \tensora \omega ) &=  ( g ( \eta \tensora - ) a^\prime \tensora g ( \omega \tensora - )  ) ( 1 - \sigma ) \nabla_0 ( \theta  )\\
&=  ( g ( \eta \tensora - ) \tensora a^\prime  g ( \omega \tensora - )  ) ( 1 - \sigma ) \nabla_0 ( \theta  )\\
&=  ( g ( \eta \tensora - ) \tensora   g ( a^\prime \omega \tensora - )  ) ( 1 - \sigma ) \nabla_0 ( \theta  ) ~ {\rm (} ~ {\rm as} ~ g ~ {\rm is} ~ {\rm bilinear} ~ {\rm )} \\
&= \psi_{a^\prime \omega, \theta} ( \eta )\\
&= \phi_\eta ( \theta \tensora a^\prime \omega ).
\end{align*}
This proves the lemma.
\qed

\bcrlre \label{21stjuly20191}
Fix $ \eta \in \Ecenter $; then there is a right $ \A $-linear map $ \widetilde{\phi}_\eta: \E \tensora \E \rightarrow \A $ defined, for all $ a \in \A, \theta, \omega \in \Ecenter $, by
$$ \widetilde{\phi}_\eta ( \theta \otimes_{\Acenter} \omega \otimes_{\Acenter} a ) = \phi_\eta ( \theta \otimes_{\Acenter} \omega ) a. $$

\noindent
Moreover, $ \widetilde{\phi}_\eta $ is a left $ \A $-linear map from $ \E \tensora \E $ to $ \A, $ i.e, $ \widetilde{\phi}_\eta \in {}_\A \Hom ( \E \tensora \E, \A ) $. Hence, by part 2. of Proposition \ref{vg2nondegenerate}, there exists a unique element $ \nabla ( \eta ) \in \E \tensora \E $ so that
\begin{equation}\label{intlambda}
\tfrac{1}{2} \widetilde{\phi}_\eta ( \xi ) = g^{(2)} (  \xi \tensora \nabla ( \eta ) ) 
\end{equation}
for all $ \xi \in \E \tensora \E $.
\ecrlre
\noindent {\bf Proof:} We only prove the statement that $ \widetilde{\phi}_\eta \in {}_\A \Hom ( \E \tensora \E, \A ) $. For this it is enough to prove that for all $ \omega, \theta \in \Ecenter $ and $ a, b  \in \A, $ we have
$$ \widetilde{\phi}_\eta ( a ( \omega \tensora \theta b ) ) = a. \widetilde{\phi}_\eta ( \omega \tensora \theta b )  .$$ 
As  $ \widetilde{\phi}_\eta $ is right $ \A $-linear by construction, it follows that
\begin{align*}
 \widetilde{\phi}_\eta ( a ( \omega \tensora \theta b ) ) &=  \widetilde{\phi}_\eta ( a ( \omega \tensora \theta ) b  ) =  \widetilde{\phi}_\eta ( a ( \omega \tensora \theta  ) ) b =  \widetilde{\phi}_\eta (  \omega \tensora \theta  a ) b \\
 &= \widetilde{\phi}_\eta (   \omega \tensora \theta  ) a b = a.  \widetilde{\phi}_\eta (   \omega \tensora \theta  )  b = a.  \widetilde{\phi}_\eta (   \omega \tensora \theta b ).
\end{align*}
In the above we have used that $ \widetilde{\phi}_\eta (   \omega \tensora \theta  ) \in \Acenter $ (
Lemma \ref{21stjuly20192} ) and $  \omega \tensora \theta ~ \in ~ \mathcal{Z} ( \E \tensora \E ) .$ 
\qed

\blmma \label{17thjuly20193}
For $ \omega, \theta, \eta \in \Ecenter, a^\prime \in \Acenter$ one has, 
\begin{equation}  \psi_{\omega, \theta} ( \eta a^\prime ) = \psi_{\omega, \theta} ( \eta ) a^\prime + 2 g ( \omega \tensora \eta ) g ( \theta \tensora d a^\prime ).  \end{equation}
\elmma
\noindent
{\bf Proof:} The proof of this Lemma  follows by a computation using the facts ( from Lemma \ref{lemma2} ) that $ da^\prime  \in \Acenter, g ( \alpha \tensora \beta ) \in \Acenter  $ for all $ \alpha, \beta \in \Ecenter $. Moreover, we also use the statement proved in Lemma \ref{16thjuly20192} that $ ( 1 - \sigma ) \nabla_0 ( e ) \in \mathcal{Z} ( \E \tensora \E ) $ for all $ e \in \Ecenter.$
We compute
\begin{align*}
&\psi_{\omega, \theta} ( \eta a^\prime )\\
&=  g ( \omega \tensora d g ( \eta a^\prime \tensora \theta ) ) - g ( \eta a^\prime \tensora d g ( \theta \tensora \omega ) ) + g ( \theta \tensora d g ( \omega \tensora \eta ) ) a^\prime \\ 
& \quad + g ( \theta \tensora g ( \omega \tensora \eta ) d a^\prime ) - ( g ( \eta \tensora - ) \tensora g ( \theta \tensora - ) ) ( ( 1 - \sigma ) \nabla_0 ( \omega )  ) a^\prime \\
& \quad+ ( g ( \omega \tensora - ) \tensora g ( \theta \tensora - ) ) ( ( 1 - \sigma ) \nabla_0 ( \eta )  ) a^\prime\\
 & \quad + ( g ( \omega \tensora - ) \tensora g ( \theta \tensora - ) ) ( ( 1 - \sigma ) ( \eta \tensora d a^\prime  ) ) - ( g ( \eta \tensora - ) \tensora g ( \omega \tensora - ) )  \nabla_0 ( \theta ) a^\prime\\
&= g ( \omega \tensora dg ( \eta \tensora \theta ) ) a^\prime + g ( \omega \tensora g ( \eta \tensora \theta ) d a^\prime    ) - g ( \eta \tensora d g (  \theta \tensora \omega  ) ) a^\prime  \\
& \quad + g ( \theta \tensora d g ( \omega \tensora \eta  ) ) a^\prime
- ( g ( \eta \tensora - ) \tensora g ( \theta \tensora - ) ) ( ( 1 - \sigma ) \nabla_0 ( \omega )  ) a^\prime 
\\ &\quad + ( g ( \omega \tensora - ) \tensora g ( \theta \tensora - ) ) ( ( 1 - \sigma ) \nabla_0 ( \eta )  ) a^\prime
- ( g ( \eta \tensora - ) \tensora g ( \omega \tensora - ) )  \nabla_0 ( \theta ) a^\prime
\\ &\quad + g ( \omega \tensora \eta ) g ( \theta \tensora d a^\prime ) - g ( \omega \tensora d a^\prime ) g ( \theta \tensora \eta )
+ g ( \omega \tensora \eta ) g ( \theta \tensora d a^\prime )\\
& = \psi_{\omega, \theta} ( \eta ) a^\prime + 2 g ( \omega \tensora \eta ) g ( \theta \tensora d a^\prime ). 
\end{align*}
This proves the lemma.
\qed

\blmma \label{20thjuly20193}
For $ e \in \Ecenter, $ let  $ \nabla ( e )   $ be the element in $ \E \tensora \E $ as defined in Corollary \ref{21stjuly20191}. Then for all $ a \in \Acenter, $ we have
\begin{equation} \label{21stjuly20193} \nabla ( e a ) = \nabla ( e ) a + e \tensora da. \end{equation}
\elmma
\noindent {\bf Proof:} It suffices to prove that for all $ \omega, \theta \in \Ecenter, $ we have
\begin{equation} \label{20thjuly20192} g^{(2)} ( ( \theta \tensora \omega ) \tensora ( \nabla ( e a ) - \nabla ( e ) a - e \tensora da   )  ) = 0. \end{equation}
Thus, since $ \{ \omega \tensora \theta: \omega, \theta \in \Ecenter  \} $ is left $ \A $-total in $ \E \tensora \E $ and $ g^{(2)} $ is left $ \A $-linear ( Proposition \ref{vg2nondegenerate} ), formula \eqref{21stjuly20193} will imply 
$$ g^{(2)} ( \xi \tensora ( \nabla ( e a ) - \nabla ( e ) a - e \tensora da   )  ) = 0 $$
for all $ \xi \in \E \tensora \E $ and hence Proposition \ref{vg2nondegenerate} will imply that
$$ \nabla ( e a ) - \nabla ( e ) a - e \tensora da  = 0.$$
For $ a \in \Acenter, $ have
$$ g^{(2)} ( ( \theta \tensora \omega ) \tensora \nabla ( e a ) ) = \tfrac{1}{2} \psi_{\omega, \theta} ( e a ).  $$
Therefore,
\begin{align*}
 g^{(2)} ( ( \theta \tensora \omega ) \tensora ( \nabla ( e a ) & - \nabla ( e ) a - e \tensora da ) ) \\
&= \tfrac{1}{2} \psi_{\omega, \theta}  ( e a ) - \tfrac{1}{2} \psi_{\omega, \theta} ( e ) a - g^{(2)} ( ( \theta \tensora \omega ) \tensora ( e \tensora da ) )\\
&=  \tfrac{1}{2} \psi_{\omega, \theta}  ( e  ) a + g ( \omega \tensora e ) g ( \theta \tensora da ) \\
& \quad - \tfrac{1}{2} \psi_{\omega, \theta} ( e ) a - g ( \theta \tensora da ) g ( \omega \tensora e ) \\
&=  0
\end{align*}
where we have used that $ g ( \omega \tensora e ) \in \Acenter $ ( Lemma \ref{lemma2} ) and Lemma \ref{17thjuly20193}.
	\qed
	
	\bppsn \label{22ndjuly20194} 
	Given the map $\nabla$ defined implicitly by formula \eqref{intlambda} and the connection $\nabla_0$ in Theorem~\ref{torsionless}, the map
	$$ L := \nabla - \nabla_0: \Ecenter \rightarrow \E \tensora \E  .$$
is right $ \Acenter $-linear and so $ L $ extends to  a right $ \A $-linear map
	$$ \widetilde{L}: \E = \Ecenter \otimes_{\Acenter} \A \rightarrow \E \tensora \E; ~ \omega \tensora a \mapsto L ( \omega ) a  .$$
	\eppsn
	\noindent {\bf Proof:} For $ \omega \in \Ecenter $ and $ a \in \Acenter, $ we have
	$$ L ( \omega a ) = ( \nabla - \nabla_0 ) ( \omega a ) = \nabla ( \omega ) a + \omega \tensora da - \nabla_0 ( \omega ) a - \omega \tensora da = L ( \omega ) a $$
	since $ \nabla_0 $ is a connection on $ \E $ and we have used Proposition \ref{20thjuly20193}. 
	\qed
	
	\bcrlre \label{21stjuly20194}
	Consider the map 
	$ \nabla: \E \rightarrow \E \tensora \E $ given 
	by the formula
	$$ \nabla = \nabla_0 + \widetilde{L} $$
	where $ \widetilde{L} $ is the map in Proposition \ref{22ndjuly20194}.
	Then $ \nabla $ is a connection on $\E$ which extends the map $ \nabla $ given in Corollary \ref{21stjuly20191}.
	\ecrlre
	\noindent {\bf Proof:} The map $ \nabla $ is a connection being the sum of a connection $ \nabla_0 $ and a right $ \A $-linear map $ \widetilde{L} $.
	\qed

\medskip
We are finally ready for the

\noindent
{\bf Proof of Theorem \ref{existenceuniqueness}} \\
The uniqueness follows from Theorem \ref{koszul21stjuly}
and the Definition of $\nabla$ in Corollary \ref{21stjuly20194}. So we are left with proving the existence. 

 We start by proving that the connection $ \nabla $ defined in Corollary \ref{21stjuly20194} is torsionless. Let   $ \widetilde{L} = \nabla - \nabla_0 $ as in Corollary \ref{21stjuly20194}. Then it suffices to prove that $ \wedge L = 0 $ since this implies
$$ \wedge \nabla = \wedge \nabla_0 = - d $$
as $ \nabla_0 $ is torsionless.
  By right $ \A $-linearity of $ \widetilde{L}, $ it suffices to check that for all $ \eta \in \Ecenter, $
	$$ \wedge \widetilde{L} ( \eta ) = 0  .$$
	However,
$$ \wedge ( \widetilde{L} ) ( \eta ) = \wedge ( \etazero \tensora \etaone - \zeroeta \tensora \oneeta ) $$
and so we need to prove that 
\begin{align*} 
\sigma ( \etazero \tensora \etaone -\zeroeta \tensora \oneeta ) &= \etazero \tensora \etaone -\zeroeta \tensora \oneeta, \\
\mbox{that is,} \qquad \etaone \tensora \etazero - \oneeta \tensora \zeroeta &= \etazero \tensora \etaone -\zeroeta \tensora \oneeta.  \end{align*}  
By bilinearity of $ g^{(2)}, $ the fact that $ \{ \theta \tensora \omega ; \theta, \omega \in \Ecenter \} $ is left $ \A$-total in $ \E \tensora \E $ and Proposition \ref{vg2nondegenerate}, it suffices to prove that for all $ \omega, \theta \in \Ecenter, $ the  following equation holds: 
$$
g^{(2)} ( ( \theta \tensora \omega  ) \tensora ( \etaone \tensora \etazero - \oneeta \tensora \zeroeta   )   ) = g^{(2)} ( ( \theta \tensora \omega  ) \tensora ( \etazero \tensora \etaone -\zeroeta \tensora \oneeta  )   ). 
$$
By a simple computation using the facts that $ g ( \omega \tensora \etaone ), g ( \theta \tensora \thetaone ) \in \Acenter,  $ 
the previous expression is seen to be equivalent to:
\begin{multline}
g ( \theta \tensora \etazero ) g ( \omega \tensora \etaone  ) - g ( \omega \tensora \etazero ) g ( \theta \tensora \etaone ) \\ = g ( \theta \tensora \zeroeta ) g ( \omega \tensora \oneeta  ) - g ( \omega \tensora \zeroeta ) g ( \theta \tensora \oneeta ).  
\end{multline}
Now, using the expressions for $ g ( \theta \tensora \etazero ) g ( \omega \tensora \etaone  ) = \tfrac{1}{2} \psi_{\theta, \omega} ( \eta ) $ and $  g ( \omega \tensora \etazero ) g ( \theta \tensora \etaone ) = \tfrac{1}{2} \psi_{\omega, \theta} $ 
( see \eqref{21stjuly20197} ) and using the facts that $ g ( \theta \tensora \oneeta ), g ( \omega \tensora \oneeta ) \in \Acenter  $ and $ g \sigma = g, $ the left hand side of the previous expression reduces to the right hand side by a straightforward simplification. 

Next we prove that $ \nabla $ is compatible with $g.$  We claim that for all $ \omega, \eta \in \Ecenter, $ we have
$$ g ( \omegazero \tensora \eta ) \omegaone + g ( \omega \tensora \etazero ) \etaone = d ( g ( \omega \tensora \eta )  ). 
$$
By virtue of Remark \ref{21stjuly201910}, this is equivalent to having for all $ \theta \in \Ecenter, $
$$
g ( \eta \tensora \omegazero ) g ( \theta \tensora \omegaone ) + g ( \omega \tensora \etazero ) g ( \theta \tensora \etaone ) = g ( \theta\tensora dg ( \omega \tensora \eta ) ) 
$$
having used the facts that $ g ( \theta \tensora - ) $ is left $ \A $-linear.  
Using the definition of $ \nabla, $ we have
\begin{align*} 
 2 g ( \eta \tensora \omegazero ) & g ( \theta \tensora \omegaone ) + 2 g ( \omega \tensora \etazero ) g ( \theta \tensora \etaone )\\
&=  g ( \eta \tensora dg ( \omega \tensora \theta ) ) - g ( \omega \tensora dg ( \theta \tensora \eta ) ) + g ( \theta \tensora dg ( \eta \tensora \omega ) )\\
 & \quad -  g ( \omega \tensora \zeroeta ) g ( \theta \tensora \oneeta ) + g  ( \omega \tensora \oneeta ) g ( \theta \tensora \zeroeta ) \\
 & \quad + g ( \eta \tensora \zeroomega ) g ( \theta \tensora \oneomega ) - g ( \eta \tensora \oneomega ) g ( \theta \tensora \zeroomega ) \\
 & \quad - g ( \omega \tensora \zerotheta ) g ( \eta \tensora \onetheta ) + g ( \omega \tensora \onetheta ) g ( \eta \tensora \zerotheta )\\
& \quad +   g ( \omega \tensora dg ( \eta \tensora \theta ) ) - g ( \eta \tensora dg ( \theta \tensora \omega ) ) + g ( \theta \tensora dg ( \omega \tensora \eta ) )   \\
& \quad - g ( \eta \tensora \zeroomega ) g ( \theta \tensora \oneomega ) + g ( \eta \tensora \oneomega ) g ( \theta \tensora \zeroomega ) \\
 & \quad + g ( \omega \tensora \zeroeta ) g ( \theta \tensora \oneeta ) - g ( \omega \tensora \oneeta ) g ( \theta \tensora \zeroeta ) \\
 & \quad - g ( \eta \tensora \zerotheta ) g ( \omega \tensora \onetheta ) + g ( \eta \tensora \onetheta ) g ( \omega \tensora \zerotheta )\\
 &= 2 g ( \theta \tensora dg ( \omega \tensora \eta ),
\end{align*}
using $ g ( \omega \tensora \onetheta ), g ( \eta \tensora \onetheta ) \in \Acenter $ and $ g ( \alpha \tensora \beta ) = g ( \beta \tensora \alpha ) $ for all $ \alpha, \beta \in \Ecenter $.
Therefore, $ \nabla $ is compatible with $ g $ on $\Ecenter$ ( as in Definition \ref{compatibilitycenter} ). This completes the proof.
 \qed

\medskip
We finish this section by comparing this result with that in \cite{article1}. For this it will be useful to adopt the notation: 
\begin{align*}
& \Pi_g^0(\nabla):\Ecenter\tensorc\Ecenter \rightarrow \E,  \\
& \Pi_g^0(\nabla)(\omega \tensorc \eta) = (g \tensora \id )\sigma_{23}(\nabla(\omega)\tensora \eta ) + ( g \tensora \id    ) ( \omega \tensora \nabla ( \eta ) ).
\end{align*}
Thus, Definition \ref{compatibilitycenter} can be rephrased by saying that $ \nabla $ is compatible with $ \E $ on $ \Ecenter $ if 
$$ \Pi_g^0 ( \nabla ) ( \omega \otimes_{\Acenter} \eta ) = d ( g ( \omega \tensora \eta ) ) ~ \forall ~ \omega, \eta \in \Ecenter  .$$
In Subsection 4.1 of \cite{article1}, it was shown that the assumption $ \E = \Ecenter \otimes_{\Acenter} \A  $ allows one to define a canonical extension $ \Pi_g ( \nabla ) : \E \tensora \E \rightarrow \E $ of the map $ \Pi^0_g ( \nabla ) $. More precisely, for all $ \omega, \eta \in \Ecenter $ and $ a \in \A, $  one has
\begin{equation} \label{22ndjuly20195} \Pi_g ( \nabla ) ( \omega \tensora \eta a ) = \Pi^0_g ( \nabla )  ( \omega \otimes_{\Acenter} \eta ) a + g ( \omega \tensora \eta ) da. \end{equation}
 We say that a connection $ \nabla $ is compatible with $ g $ on the whole of $ \E $ if for all $e, f$ in $\E,$
 \begin{equation} \label{10thaugust20191} \Pi_g ( \nabla ) ( e \tensora f ) = d ( g ( e \tensora f ) ). \end{equation} 
It was also shown that for any connection $ \nabla_1 $ on $ \E, $   the map from $ \E \tensora \E $ to $ \E,$ defined by
$$ ~ e \tensora f \mapsto \Pi_g ( \nabla_1 ) ( e \tensora f ) - d ( g ( e \tensora f ) ) $$
is right $ \A $-linear.
We have the following result which recovers the main result of \cite{article1}:

\bcrlre
Suppose $ \E:= \oneform $ satisfies the assumptions of Theorem \ref{existenceuniqueness} and $ g $ is a pseudo-Riemannian bilinear metric on $ \E $. Then there exists a unique connection on $ \E $ which is torsionless and compatible with $ g $ on the whole of $\E.$
\ecrlre  
  \noindent {\bf Proof:} From Theorem \ref{existenceuniqueness}, we know that there exists a unique connection $ \nabla $ which is torsionless and compatible with $  g$ on $ \Ecenter.$ Thus, for all $ \omega, \eta \in \Ecenter, $
	$$ \Pi^0_g ( \nabla ) ( \omega \otimes_{\Acenter} \eta ) = d ( g ( \omega \tensora \eta ) )  .$$
	Therefore, \eqref{22ndjuly20195} implies that for all $ a \in \A, $
	$$ \Pi_g ( \nabla ) ( \omega \tensora \eta a ) = d ( g ( \omega \tensora \eta )  ) a + g ( \omega \tensora \eta ) da = d ( g ( \omega \tensora \eta a ) ), $$
	i.e, $ \nabla $ is compatible with $ g $ on the whole of $ \E $. Uniqueness is clear from Theorem \ref{existenceuniqueness}.
	\qed

\section{Levi-Civita connections as bimodule connections} \label{section7}

In this section, we make contact with bimodule connections. 
A considerable amount of literature
on Levi-Civita connections in the context of noncommutative geometry have been devoted to bimodule connections. We refer to the book \cite{beggsmajidbook} for the details. We show that in our set up ( the assumptions of Theorem \ref{existenceuniqueness} ), the Levi-Civita connection is a bimodule connection in a very natural way. This section is a genuine application of the Koszul formula of Theorem  \ref{koszul21stjuly}. Let us recall the definition of bimodule connections.  

\bdfn\label{lLr}
Suppose $ \E = \oneform $ for a differential calculus and $ \sigma^\prime: \E \tensora \E \rightarrow \E \tensora \E  $ be a bimodule map. A right connection $\nabla_1$ on $  \E $ is said to be a bimodule connection for the pair $ ( \E, \sigma^\prime ) $ if, in addition to the right Leibniz rule as in Definition \ref{rLr},  there is also a $\sigma^\prime$-left Leibniz rule, that is,  for all $ a \in \A $ and for all $ e \in \E, $ it holds that
$$ \nabla_1 ( a e  ) = a \nabla_1 ( e ) + \sigma^\prime ( da \tensora e )  .$$  
\edfn
Throughout this section, we will work under the assumptions of Theorem \ref{existenceuniqueness} and so we have a canonical choice of $ \sigma^\prime = \sigma $ as defined in Definition \ref{16thjuly201923}. We will show that the Levi-Civita connection of Theorem \ref{existenceuniqueness} is a bimodule connection for the pair $ ( \E, \sigma ) $.

We start by proving a necessary and sufficient condition for a connection to be a bimodule connection for the pair $ ( \E, \sigma ) $.  

\bppsn \label{16thjuly20193}
Suppose $ \E = \oneform $ as in Theorem \ref{existenceuniqueness}. If $ \nabla_1 $ is a connection on $ \E, $ then it is a bimodule connection for $ ( \E, \sigma ) $ if and only if $ \nabla_1 ( \Ecenter ) \subseteq \mathcal{Z} ( \E \tensora \E ).$
\eppsn
\noindent {\bf Proof:} Suppose $  \nabla_1 $ is a bimodule connection for the pair $ ( \E, \sigma ) $. Then for all $ \omega \in \Ecenter $ and $ a \in \A, $ we get
\begin{align*}
 \nabla_1 ( \omega ) a + \omega \tensora da &= \nabla_1 ( \omega a ) = \nabla_1 ( a \omega )\\
&= a \nabla_1 ( \omega ) + \sigma ( da \tensora \omega ) = a \nabla_1 ( \omega ) + \omega \tensora da,  
\end{align*}
since $ \omega \in \Ecenter $. Thus, for all $ \omega \in \Ecenter $ and $ a \in \A, $ we have
$$ \nabla_1 ( \omega ) a = a \nabla_1 ( \omega ), $$
that is, $ \nabla_1 ( \Ecenter ) \subseteq \mathcal{Z} ( \E \tensora \E ) $.

Conversely, suppose $ \nabla_1 $ is a connection such that $ \nabla_1 ( \Ecenter ) \subseteq \mathcal{Z} ( \E \tensora \E ) $. Then for all $ \omega \in \Ecenter $ and for all $ a \in \A, $
\begin{equation} \label{12thaugust20191} \nabla_1 ( a \omega ) = \nabla_1 ( \omega a ) = \nabla_1 ( \omega ) a + \omega \tensora da = a \nabla_1 ( \omega ) + \sigma ( da \tensora \omega ).\end{equation} 
Now let $ e \in \E $ and $ a \in \A $. Since $ \E $ is centered, we can write $ e = \sum\nolimits_j f_j b_j $ for some $ f_j \in \Ecenter $ and $ b_j \in \A $. Then
\begin{align*}
\nabla_1 ( a e ) &= \sum\nolimits_j \nabla_1 ( a f_j b_j ) = \sum\nolimits_j [ \nabla_1 ( a f_j ) b_j + a f_j \tensora d b_j ]\\
&= \sum\nolimits_j [ a \nabla_1 ( f_j ) b_j + \sigma ( da \tensora f_j ) b_j + a f_j \tensora d b_j ] \qquad {\rm (} ~ {\rm by} ~ \eqref{12thaugust20191} ~ {\rm )}\\
 &= \sum\nolimits_j [ a \nabla_1 ( f_j b_j ) + \sigma ( da \tensora f_j ) b_j ]\\
&=  a \nabla_1 ( e ) + \sigma ( da \tensora e ).
\end{align*}  
Therefore, $ \nabla_1 $ is a bimodule connection.
\qed

Now we use the Koszul formula to prove the main result of this subsection. 

\bthm \label{16thjuly20194}
Suppose $ \E = \oneform $ satisfies the assumptions of Theorem \ref{existenceuniqueness} and $ g $ is a pseudo-Riemannian bilinear metric on $\E$. Then the Levi-Civita connection for $ ( \E, g ) $ obtained in Theorem \ref{existenceuniqueness} is a bimodule connection for the pair $ ( \E, \sigma ) $.
\ethm
\noindent {\bf Proof:} The proof follows from Proposition \ref{16thjuly20193} and the Koszul formula as expressed in equation \eqref{16thjuly2019koszul}. We claim that it is enough to show that for all $ \omega, \eta, \theta \in \Ecenter, $
\begin{equation} \label{23rdjuly20191} 2 ( g ( \omega \tensora -  ) \tensora  g ( \theta \tensora -  )  ) ( \nabla ( \eta ) ) \in \Acenter. \end{equation}
Indeed, by virtue of Remark \ref{11thaugust2019}, we have
\begin{align*}
( g ( \omega \tensora -  ) & \tensora g ( \theta \tensora -  )  ) ( \nabla ( \eta ). a - a. \nabla ( \eta )   ) \\ & = g^{(2)} ( ( \theta \tensora \omega ) \tensora ( \nabla ( \eta ). a - a. \nabla ( \eta )   ) )\\
 &= g^{(2)} ( ( \theta \tensora \omega ) \tensora \nabla ( \eta ) ). a - g^{(2)} ( ( \theta \tensora \omega ). a \tensora \nabla ( \eta)  )\\
&= g^{(2)} ( ( \theta \tensora \omega ) \tensora \nabla ( \eta ) ). a - g^{(2)} ( a. ( \theta \tensora \omega ) \tensora \nabla ( \eta)  )\\
&= g^{(2)} ( ( \theta \tensora \omega ) \tensora \nabla ( \eta ) ). a - a. g^{(2)} ( ( \theta \tensora \omega ) \tensora \nabla ( \eta)  )\\
&= ( g ( \omega \tensora -  ) \tensora g ( \theta \tensora -  )  ) ( \nabla ( \eta ) ). a - a. ( g ( \omega \tensora -  ) \tensora g ( \theta \tensora -  )  ) ( \nabla ( \eta ) ),
\end{align*}
which is equal to zero when \eqref{23rdjuly20191} holds. Since this is true for all $ \theta, \omega \in \Ecenter, $ a combination of part 4. of  Lemma \ref{20thaugust20191}, Remark \ref{11thaugust2019}  and the left $\A$-linearity imply that 
$$ g^{(2)} ( \xi \tensora ( \nabla ( \eta ). a - a. \nabla ( \eta )   ) ) = 0 ~ \forall \xi \in \E \tensora \E  .$$  
  This proves our claim by Lemma \ref{vg2nondegenerate}.
	
	Now, by Lemma \ref{lemma2}, we have that the elements $ g ( \omega \tensora  dg ( \eta \tensora \theta ) )$, $g ( \eta \tensora dg ( \theta \tensora \omega  ) )$ and $g ( \theta \tensora dg ( \omega \tensora \eta ) )$ are all in $\Acenter $. Next, we observe that for all $ a $ in $ \A, $
	\begin{align*}
	 ( g ( \eta \tensora -  ) \tensora g ( \theta \tensora - ) ) ( 1 - \sigma ) \nabla_0 ( \omega ). a  &= g^{(2)} ( ( \theta \tensora \eta ) \tensora ( ( 1 - \sigma ) \nabla_0 ( \omega ). a ) )\\
	&= g^{(2)} (  ( \theta \tensora \eta ) \tensora ( a. ( 1 - \sigma ) \nabla_0 ( \omega ) ) )\\
	&= a. g^{(2)} ( ( \theta \tensora \eta ) \tensora ( ( 1 - \sigma ) \nabla_0 ( \omega ) ) )\\
	&= a. ( g ( \eta \tensora -  ) \tensora g ( \theta \tensora -  ) ) ( 1 - \sigma ) \nabla_0 ( \omega ).
	\end{align*}
	Here we have used the fact that $ ( 1 - \sigma ) \nabla_0 ( \omega ) \in \mathcal{Z} ( \E \tensora \E ) $ by   Lemma \ref{16thjuly20192}. This proves that $ ( g ( \eta\tensora -  ) \tensora g ( \theta \tensora -  ) ) ( 1 - \sigma ) \nabla_0 ( \omega ) \in \Acenter $. Similarly, one finds that also $ (  g ( \omega \tensora -  ) \tensora g ( \theta \tensora -  ) ) ( 1 - \sigma ) \nabla_0 ( \eta ) $ and $  ( g ( \eta \tensora -  ) \tensora g ( \omega \tensora -  )   ) ( 1 - \sigma ) \nabla_0 ( \theta )$ belong to $\Acenter$. Hence, by the Koszul formula given in \eqref{16thjuly2019koszul}, $ ( g ( \omega \tensora -  ) \tensora g ( \theta \tensora -   ) ) \nabla ( \eta ) \in \Acenter $ which completes the proof of the the theorem. 
\qed

\section{The example of the fuzzy  sphere} \label{section8}
This section concerns the example of a spectral triple on the fuzzy  sphere. Our goal is to prove the existence and uniqueness of Levi-Civita connections ( on the corresponding module of one-forms ) for any bilinear pseudo-Riemannian metric. Our spectral triple is a truncated version of the spectral triple constructed in \cite{frolich} for a  fuzzy $3$-sphere.

It turns out that the module of one-forms is free of rank $3.$ We compute the connection forms (and the Christoffel symbols) of the Levi-Civita connection for the canonical pseudo-Riemannian metric coming from the spectral triple ( see \cite{frolich} ). The computations in this section are similar to those in Section 3 of \cite{frolich} and those in Section 5 of \cite{article1}. However, we provide all the details for the sake of completeness.
  
Let us set up some notations. Firstly, the Lie algebra of $so(3) \simeq su(2)$ is generated by three elements $J_k, k=1,2,3$ with commutation relations
\begin{equation}\label{23rdjuly20193}
[ J_k, J_l ] = \sum\nolimits_{m=1}^3 \epsilon_{klm} J_m. 
\end{equation}
Here $\epsilon_{klm} $ is the completely antisymmetric Levi-Civita rank $3$ tensor with $\epsilon_{123} = 1$.

For a natural number $n,$ let $ \rho_{\frac{n}{2}} $ denote the $(n + 1)$-dimensional unitary irreducible representation of the Lie algebra $su(2)$. The vector space $ \mathbb{C}^{n + 1} $ is the carrier vector space of the representation $ \rho_{\frac{n}{2}}$, and $ \mathbb{C} $ will be viewed as the trivial representation space. In particular, for the fundamental, $n=1$, representation we have 
$$
J_k = \tfrac{1}{2} \sqrt{-1} \, \tau_k, 
$$
with Hermitian Pauli matrices 
$$
\tau_1 = \begin{bmatrix} 0 & 1 \\ 1 & 0  \end{bmatrix}, 
\quad \tau_2 = \begin{bmatrix} 0 & - \sqrt{-1} \\ \sqrt{-1} & 0  \end{bmatrix}, 
\quad \tau_3 = \begin{bmatrix} 1 & 0 \\ 0  & -1 \end{bmatrix},
$$
which are also a basis of the Clifford algebra $Cl(2,0)$, that is are such that,
\begin{equation} \label{23rdjuly20192} 
[\tau_k, \tau_l ] = 2 \sqrt{-1}\, \sum\nolimits_{m=1}^3 \epsilon_{klm} \tau_m , \qquad \tau_j \tau_k + \tau_k \tau_j = 2 \delta_{j k} .
\end{equation}

Our spectral triple $ ( A_N, H_N, D_N ) $ has Hilbert space 
$$ 
H_N = K_N \tensorc \mathbb{C}^2, 
$$ 
where 
$ K_N =  \oplus_{l=0}^{N} \mathbb{C}^{2l+1}$. 
The algebra $ A_N $ is the full matrix algebra $ B ( K_N ) $.
We have the canonical  action $ \pi^\prime $ of $ A_N $ on $ K_N $. Then  
the algebra  $ A_N $ is represented on $ H_N  $ by the formula $ a \mapsto \pi ( a ) $ where 
$$ \pi ( a ) ( h_1 \tensorc h_2 ) = \pi^\prime (a h_1) \tensorc h_2 $$
for all $ h_1 $ in $ K_N $ and $ h_2 $ in $ \mathbb{C}^2.$ 

Next we have the Dirac operator $ D_N $.  
Firstly, for $k = 1,2,3, $ we define operators $ X_k \in B ( H_N ) $ by the formula
 $$ 
 X_k = \oplus^{N}_{n = 0} \,\, \rho_{\frac{n}{2}} ( J_k ),
 $$
 and denote $\sigma_k = \sqrt{-1} \, \tau_k$,  for $k = 1,2,3$. 
Then
the Dirac operator $ D_N $ is defined as
$$ D_N = \sum^3_{k = 1}   X_k  \tensorc \sigma_k. $$
Since $ X_k $ and and $ \sigma_k $ are all skew-Hermitian, the operator $ D_N $ is self-adjoint.

By omitting the notation $ \pi $ while viewing an element of $ A_N $ as an operator on $ H_N $, it can be easily checked that for all $ a \in A_N, $
\begin{equation} \label{23rdjuly20194} [ D_N, a ] = \sum^3_{k = 1} \, [ X_k, a ] \tensorc \sigma_k, \end{equation}
as elements of 
$ A_N \tensorc M_{2} ( \mathbb{C})$ acting on $H_N$.
For $ k = 1,2, 3, $ we define derivations $ \delta_k $ on $ A_N $ by 
  $\delta_k ( a ) = [ X_k, a ]$.
	Then
	\begin{equation} \label{28thjuly20194} [ D, a ] = \sum^3_{k = 1} \, \delta_k ( a ) \tensorc \sigma_k. \end{equation}
	By use of \eqref{23rdjuly20193}, one has the following commutation relations between the derivations $\delta_i $:
	\begin{equation} \label{28thjuly20197} [ \delta_1, \delta_2 ] = \delta_3 , \quad [ \delta_2, \delta_3 ] = \delta_1, \quad 
	[ \delta_3, \delta_1 ] = \delta_2. \end{equation}
	
	We will denote the space of one forms of this spectral triple by the symbol $ \E $. Our goal is to  prove that $ \E $ satisfies all the conditions of Theorem \ref{existenceuniqueness}. We will repeatedly use the fact that 
	$A_{N} = B (K_N)$ has no proper ideal except $ \{ 0 \} $ and itself. Let us  start by identifying the space of one forms $\E$ as well as the space of two-forms. 
	We have the following proposition:
	\bppsn \label{28thjuly20192}
	The module $\E = \Omega^1(A_N) $ is free of rank $ 3 $ generated by the central elements $ 1 \tensorc \sigma_1, 1 \tensorc \sigma_2, 1 \tensorc \sigma_3 $.  
	\eppsn
	\noindent {\bf Proof:} We use \eqref{23rdjuly20193} and \eqref{23rdjuly20194} to write:
	\begin{align*} [ D, J_1 ] & = J_2 \tensorc \sigma_3 - J_3 \tensorc \sigma_2 , \\ [ D, J_2 ] & = J_3 \tensorc \sigma_1 - J_1 \tensorc \sigma_3, \\ [ D, J_3 ] & = J_1 \tensorc \sigma_2 - J_2 \tensorc \sigma_1 . 
	\end{align*}
	Hence, we get
	\begin{equation} \label{23rdjuly20195} [ J_1, [ D, J_1 ] ] = J_2 \tensorc \sigma_2 + J_3 \tensorc \sigma_3, \end{equation}
	\begin{equation} \label{23rdjuly20196} [ J_2, [ D, J_2 ] ] = J_1 \tensorc \sigma_1 + J_3 \tensorc \sigma_3, \end{equation}
	\begin{equation} \label{23rdjuly20197} [ J_3, [ D, J_3 ] ] = J_1 \tensorc \sigma_1 + J_2 \tensorc \sigma_2. \end{equation}
By ( \eqref{23rdjuly20195} - \eqref{23rdjuly20197} ) + \eqref{23rdjuly20196}, we have
$$ [ J_1, [ D, J_1 ] ] - [ J_3, [ D, J_3 ] ] + [ J_2, [ D, J_2 ] ] = 2 J_3 \tensorc \sigma_3  .$$
	Therefore, the element $J_3 \tensorc \sigma_3 $ is an element of $ \E .$ Since $ \E $ is a bimodule and the ideal generated by the non-zero element  $ J_3 $ is equal to $ A_N, $ we can conclude that $ 1 \tensorc \sigma_3 $ belongs to $ \E $. Similarly, the elements $ 1 \tensorc \sigma_1 $ and $ 1 \tensorc \sigma_2 $ also belong to $ \E $. 
	Thus, $ \E \simeq A_N \tensorc \mathbb{C}^3 $ is a free module with a basis consisting of the central elements 
	$ 1 \tensorc \sigma_k, \, k = 1,2,3 $. 
	\qed
	
	\bcrlre \label{28thjuly20191}
	Let denote by the symbol $ \E. \E $ the subset $ \{ e. f : e, f ~ \in \E     \}  $ of $ A_N \tensorc M_2(\IC)$. 
	Then $ \E. \E = A_N \tensorc M_2(\IC) $. 
	\ecrlre

 \noindent {\bf Proof:} By Proposition \ref{28thjuly20192}, it follows that the elements $ 1 \tensorc \sigma_k $ belongs to $ \E $ for all $k = 1,2, 3.$ Thus, the elements $ 1 \tensorc \sigma_k \sigma_j \in \E. \E $.  Since \eqref{23rdjuly20192} holds and $ \{ 1, \sigma_k : k = 1,2, 3 \} $ is a basis of $ B ( \mathbb{C}^2 ) = M_2(\IC)$, we see that the subspace 
 $ \{ 1 \tensorc X : X \in M_2(\IC) \} $ belongs to $ \E. \E $. Therefore, $ A_N \tensorc M_2(\IC)$ belongs to  $ \E. \E $. Hence, $ \E. \E = A_N \tensorc M_2(\IC).$   
\qed

\bppsn \label{28thjuly20193}
 The space of junk forms {\rm(} ~ $ A_N \tensorc M_2(\IC) $ ~ {\rm )} of the spectral triple is equal to the subspace $ \{ Y \tensorc 1 : Y \in A_N  \} $.
\eppsn
\noindent {\bf Proof:} Suppose $ a_k, b_k, \, k = 1, \cdots, n $, be elements in $ A_N $ such that $ \sum\nolimits_k a_k [ D, b_k ] = 0.$ Then by \eqref{28thjuly20194}, we get
$$ \sum\nolimits_{j,k} a_j \, \delta_k ( b_j ) \tensorc \sigma_k = 0  .$$
 Therefore, for all $k = 1,2,3,$ we have
\begin{equation} \label{28thjuly20195} \sum\nolimits_j a_j \, \delta_k ( b_j ) = 0. \end{equation}
We apply $ \delta_l $ to \eqref{28thjuly20195} to obtain
$$ \sum\nolimits_j \, [ \delta_l ( a_j ) \delta_k ( b_j ) + a_j \delta_l \delta_k ( b_j )  ] = 0 \qquad \forall ~ k,l,$$   
where we have used the fact that $ \delta_l $ is a derivation. Hence, for all $k,l,$ we get
\begin{equation} \label{28thjuly20196} \sum\nolimits_j \delta_l ( a_j ) \delta_k ( b_j ) = - \sum\nolimits_j a_j \delta_l \delta_k ( b_j ), \qquad \sum\nolimits_j \delta_k ( a_j ) \delta_k ( b_j ) = - \sum\nolimits_j a_j \delta^2_k ( b_j ). \end{equation}
Now, \eqref{28thjuly20196} implies that 
\begin{align*}
\sum\nolimits_j  [ D, a_j ] [ D, b_j ] & = - \sum\nolimits_j \sum^3_{k,l = 1} \delta_k ( a_j ) \delta_l ( b_j ) \tensorc \sigma_k \sigma_l\\
& = - \sum\nolimits_j   \sum^3_{k = 1} \delta_k ( a_j ) \delta_k ( b_j ) \tensorc 1  \\ & \quad - \sum\nolimits_j \sum\nolimits_{k < l} [ \delta_k ( a_j ) \delta_l ( b_j ) - \delta_l ( a_j ) \delta_k ( b_j ) ] \tensorc \sigma_k \sigma_l\\
& =   \sum\nolimits_{i,k} a_j \delta^2_k ( b_j ) \tensorc 1 + \sum\nolimits_{j, k < l} a_j [ \delta_k, \delta_l ] ( b_j ) \tensorc \sigma_k \sigma_l
 \end{align*}
However, by virtue of \eqref{28thjuly20197}, we have 
\begin{align*}
\sum\nolimits_{i, k < l} a_j [ \delta_k, \delta_l ] ( b_j ) \tensorc \sigma_k \sigma_l
& = \sum\nolimits_j a_j \delta_3 ( b_j ) \tensorc \sigma_1 \sigma_2 \\ & \quad + \sum\nolimits_j a_j \delta_1 ( b_j ) \tensorc \sigma_2 \sigma_3 - \sum\nolimits_j a_j \delta_2 ( b_j ) \tensorc \sigma_1 \sigma_3\\
&= 0  
\end{align*}
by \eqref{28thjuly20195}.
Thus, we have
\begin{equation} \label{28thjuly20198} \sum\nolimits_j [ D, a_j ] [ D, b_j ] =  \sum\nolimits_{j,k} a_j \delta^2_k ( b_j ) \tensorc 1. \end{equation}
We claim that there exist elements $a, b \in A_N$ such that $ a \delta_k ( b ) = 0 $ for all $k,$ but $ a \sum\nolimits_k \delta^2_k ( b ) \neq 0 $. Indeed, if our claim is true, then by the arguments made above, the set of all junk forms will be of the form $ \{  X \tensorc 1: X \in \mathcal{I}  \}, $ where $ \mathcal{I} $ is a non-zero ideal of $ M_{N + 2} ( \mathbb{C} ) $. Therefore, the space of  junk forms is equal to  $\{  X \tensorc 1: X \in M_{N + 2} ( \mathbb{C} )  \} $.

So we are left to prove our claim. Let us define the bilinear form on $ \mathbb{C}^{N + 1} $
$$ \left\langle \left\langle v_1, v_2 \right\rangle\right\rangle = {\rm Re} ( \left\langle v_1, v_2 \right\rangle  ), $$
where $ \left\langle  \cdot, \cdot \right\rangle $ is the complex inner product on $ \mathbb{C}^{N + 1} $ with respect to which $ J_1, J_2, J_3 $ are skew hermitian elements of $ {\rm Hom}_{\mathbb{C}} ( \mathbb{C}^{N + 1}, \mathbb{C}^{N + 1} ).$  
 
Moreover, let $ v $ be a vector in $ \mathbb{C}^{N + 1} $ such that $   J_1  ( v ) \neq 0 $. We claim that $ v  $ does not belong to ${\rm Span} \{ J_1 v, J_2 v, J_3 v  \}.$ Indeed, since $ J_k $ is skew-hermitian, we have
 $$ \left\langle  J_k v, v \right\rangle = - \left\langle v,J_k v \right\rangle = - \overline{\left\langle J_k v, v \right\rangle} $$
and hence, $ \left\langle \left\langle  J_k v, v \right\rangle\right\rangle = 0 $. From here, it is  straightforward to verify our claim. 

We construct a basis $ \{ v, v_1, \cdots v_N \} $   of $ \mathbb{C}^{N + 1} $ in such a way that $v$ is the first element and $ {\rm Span} \{ J_1 v, J_2 v, J_3 v \} \subseteq \{ v_1, \cdots v_N \} $.  
We define 
$$ a ( \lambda. 1 + x + c_1 v + \sum^N_{i = 1} c_i v_i ) = c_1 v, ~ b ( \lambda. 1 + x + c_1 v + \sum^N_{i = 1} c_i v_i ) = \lambda. v$$
for all $ \lambda $  in $ \mathbb{C} $ and for all $ x $ in $ \mathbb{C}^2 \oplus \mathbb{C}^3 \oplus \cdots \oplus \mathbb{C}^N.$ 
Then 
\begin{align*}
a \delta_k ( b ) ( \lambda. 1 + x + c_1 v + \sum^N_{i = 1} c_i v_i ) &= a ( X_k b - b X_k  ) ( \lambda. 1 + x + c_1 v + \sum^N_{i = 1} c_i v_i )\\
&= a ( X_k. b ) ( \lambda. 1 ) = \lambda a X_k ( v ) = \lambda a J_k ( v ) = 0.
\end{align*} 
 However, 
\begin{align*}
a ( \sum\nolimits_k \delta^2_k ( b )  ) ( 1 ) &= a \sum\nolimits_k [ X_k, [ X_k, b ] ] ( 1 )\\
&= a \sum\nolimits_k [ X_k,  X_k b - b X_k   ) ] ( 1 ) = a ( \sum\nolimits_k X^2_k b - X_k b X_k - X_k b X_k + b X^2_k ) ( 1 )\\
&= a (  \sum\nolimits_k X^2_k b ) ( 1 ) = a ( \sum\nolimits_k X^2_k ) ( v )\\
&= a ( \sum\nolimits_k J^2_k  ) ( v ) = - \tfrac{3}{4} a ( v ) \neq 0.
\end{align*}
This finishes the proof of our claim about the description of the junk forms.
 \qed

\bcrlre \label{17thaugust20191}
The space of two forms is
$$ \Omega^2(A_N) =  \{ X \tensorc \sigma_1 \sigma_2 + Y \tensorc \sigma_2 \sigma_3 + Z \tensorc \sigma_1 \sigma_3: X, Y, Z \in A_N  \} .
$$
 \ecrlre
\noindent {\bf Proof:} The proof is an immediate consequence of Corollary \ref{28thjuly20191} and Proposition \ref{28thjuly20193}. \qed

Let us denote by $ e_k = 1 \tensorc \sigma_k$, for $k =1,2,3$. Then, 
From Proposition \ref{28thjuly20192}, the element $ \{ e_1, e_2, e_3 \} \subseteq \Ecenter $ 
form a basis of the free right $A_N$-module $\E.$

We are in the position to prove the main theorem of this section.

\bthm 
The differential calculus coming from the spectral triple on the fuzzy  sphere satisfies the hypothesis of Theorem \ref{existenceuniqueness} and hence there exists a unique Levi-Civita connection for any pseudo-Riemannian bilinear metric on 
$\E = \Omega^1(A_N) $. 
\ethm
\noindent {\bf Proof:} Firstly, the description of $ \E $ in Proposition \ref{28thjuly20192} implies that $ \Ecenter $ is the complex linear span of $ \{ e_1, e_2, e_3 \} $. Thus,  the equality $ \Ecenter \otimes_{\ancenter} A_N = \E $ easily follows from the description of one-forms, observing that $ \ancenter = \mathbb{C}.1 $. 

Next, Corollary \ref{17thaugust20191} implies that  $ \Omega^2(A_N) = {\rm Span}\{ a_{kj} \, e_k \wedge e_j: k \neq j, a_{kj} = - a_{jk} \} $ and $ {\rm Ker} ( \wedge ) = {\rm Span}\{ a_{kj} \, e_k \tensoran e_j: a_{kj} =  a_{jk} \} $. If $ \F = {\rm Span}\{ a_{kj} \, e_k \tensoran e_j: k \neq j,  a_{kj} = - a_{jk} \}, $ then  $ \F $ is isomorphic to $ \Omega^2(A_N) $ as right $A_N$-modules and
 $$\E \tensoran \E = {\rm Ker} ( \wedge ) \oplus \F.$$
 Finally, it is easy check that for all $\omega, \eta \in \Ecenter, $   
			$$ \Psym ( \omega \tensoran \eta ) = \tfrac{1}{2} ( \omega \tensoran \eta + \eta \tensoran \omega ) $$
			and therefore, $\sigma ( \omega \tensoran \eta ) = \eta \tensoran \omega  $ for all $ \omega, \eta \in \Ecenter $. 
This verifies all the hypotheses of Theorem \ref{existenceuniqueness}.
			\qed
	
By construction, the bimodule of 2-forms is free and we take as generators the elements
\begin{equation} \label{7thaugust20192nn} 
f_m = \tfrac{1}{2} \, \sum\nolimits_{j,k} \epsilon_{mjk} e_j \wedge e_k ,  \qquad m=1,2,3. 
\end{equation}
It is then easy to deduce that 
\begin{equation} \label{7thaugust20192n} 
e_j \wedge e_k =  \sum\nolimits_{m}\, \epsilon_{jkm} f_m ,  \qquad j,k=1,2,3. 
\end{equation}
As for the action of the differential, by definition  
	$d ( a ) = \sum\nolimits_k [J_k, a] \, e_k$, for $a\in\A_N$. On the other hand, starting from $0=\delta^2 (J_k)$, a direct computation leads to
	\begin{equation} \label{7thaugust20192} 
d e_m  =  - \tfrac{1}{2} \, \sum\nolimits_{j,k} \epsilon_{mjk} \, e_j \wedge e_k = - f_m . 
\end{equation}

\bthm\label{lc-fuzzy}
The connection $1$-forms of the Levi-Civita connection for the canonical pseudo-Riemannian bilinear metric for the spectral triple are given by 
\begin{equation} \label{16thaugust20191n}
\omega_{jk} = - \tfrac{1}{2} \sum\nolimits_l \, \epsilon_{jkl} \, e_l \, .
\end{equation} 
\ethm
\noindent {\bf Proof:} It is easy to see ( as in equation (3.49) of \cite{frolich} ) that the canonical pseudo-Riemannian bilinear metric $g$ is defined as the right $A_N$-linear extension of the map  
   $$ g ( e_k \tensoran e_j ) = \delta_{kj}  .$$
	We are going to compute the Levi-Civita connection for this $g.$
 Let $ \nabla $ be a connection on $\E$ which is both torsionless and compatible with $g.$ Since $ \E $ is a free centered module with a basis $ \{ e_1, e_2, e_3 \} \subseteq \Ecenter, $ we can write 
\begin{equation} \label{16thaugust20191} 
\nabla ( e_k ) = \sum\nolimits_{j} e_j \tensoran \omega_{jk}.
\end{equation} 
Now, since the basis $e_k$ is orthonormal, the metric compatibility condition reduces to 
$$ 
(g \tensoran {\rm id})\sigma_{23}(\nabla(e_k)\tensoran e_j ) + ( g \tensoran {\rm id} ) ( e_k \tensoran \nabla ( e_j ) ) = 0.
$$
By inserting \eqref{16thaugust20191}, this just give that for all $k, j$ we have antisymmetry
\begin{equation} \label{15thaugust20192} 
\omega_{kj} = - \omega_{jk}.
 \end{equation}
Next, the  torsion zero condition $d e_k + \wedge \circ \nabla ( e_k ) = 0$ becomes 
\begin{equation}  
d e_k + \sum\nolimits_{j} e_j \wedge \omega_{jk} = 0. 
\end{equation}
Using \eqref{7thaugust20192} and the antisymmetry \eqref{15thaugust20192} one infers that 
$\omega_{jk} = - \tfrac{1}{2} \sum\nolimits_l \, \epsilon_{jkl} \, e_l$ as stated in \eqref{16thaugust20191n}. 
This finishes the proof. \qed 

\brmrk
If we write the connection $1$-forms on the basis $e_k$, 
$$
\omega_{jk} = \sum\nolimits_{l} e_l \, \Gamma_{k j l},
$$
we get for the Christoffel symbols $\Gamma_{l k j}$ of the Levi-Civita connection the expressions 
$$ 
\Gamma_{ljk} = \tfrac{1}{2} \epsilon_{ljk}.
$$
\ermrk

\subsection{Computation of the curvature}
Let us recall the Sweedler-like notation in \eqref{sweedler} for a general torsionless connection $ \nabla $ on a 
centered $A$-bimodule 
$\E: \nabla ( \omega ) = \omegazero \tensora \omegaone$. We define the right $A_N$-linear map 
$$  
R_\nabla := \nabla^2 : \, \E \rightarrow \E \tensoran \Omega^2(A_N)  
$$
by a Leibniz rule:
\begin{align}
R_\nabla(\omega) &= \nabla (\omegazero \tensoran \omegaone) 
= \nabla (\omegazero) \wedge \omegaone + \omegazero \tensoran d \omegaone \nn \\
& = (\omegazero)_{(0)} \tensoran (\omegazero)_{(1)} \wedge \omegaone + \omegazero \tensoran d \omegaone.  
\end{align}

Then, let us specialise this to the fuzzy sphere with free bimodule of 1-forms $\E$ and basis $(e_k, k=1,2,3)$ as before.
We have the following proposition.
\bppsn \label{curv-fuzzy2n}
Consider the Levi-Civita connection of Theorem \ref{lc-fuzzy} with connections 
$1$-forms $\omega_{jk} = - \frac{1}{2} \sum\nolimits_l \, \epsilon_{jkl} \, e_l$. Then
\begin{align}\label{curv}
R_\nabla(e_j) 
&=
- \tfrac{1}{4} \sum\nolimits_{p,q} \, \varepsilon_{jpq} \ e_p \tensoran f_q  \nn \\
&=
\tfrac{1}{4} \sum\nolimits_{p} \, e_p \tensoran ( e_p \wedge e_j ) \, .
\end{align}
with basis 1-forms $e_p$'s and basis two forms $f_q$'s.
\eppsn
\noindent {\bf Proof:} 
Writing as before $\nabla ( e_k ) = \sum\nolimits_{j} e_j \tensoran \omega_{jk}$ for the connection,
one finds
$$
R_\nabla (e_j) = \sum\nolimits_{p,k} \, e_p \tensoran \left( \omega_{p k} \wedge \omega_{k j}  + d ( \omega_{p j} ) \right) 
$$
Then, a direct computation using \eqref{16thaugust20191n}, gives
$$
\omega_{p k} \wedge \omega_{k j} = \tfrac{1}{4} \sum\nolimits_{q} \, \varepsilon_{jpq} \, f_q \quad 
\mbox{and} \quad d \omega_{p j} = - \tfrac{1}{2} \sum\nolimits_{q} \, \varepsilon_{jpq} \, f_q
$$ 
which when inserted in the previous expression lead to the first expression in \eqref{curv}.
Using equation \eqref{7thaugust20192nn} one gets the second equality.
\qed

Since the bimodule $\E$ is free with basis of generators $(e_k, k=1,2,3)$ so is the dual bimodule $\E^*$ whose dual basis we shall denote $(\psi_j, j=1,2,3)$ with $\psi_j(e_k ) = \delta_{jk}$. Then, by using the canonical isomorphism 
$$
\Hom_{A_N}(\E, \E \tensoran \Omega^2(A_N)) \simeq \E \tensoran \Omega^2(A_N) \tensoran \E^*
$$
we can think of the curvature $R_\nabla$ as a map 
$R_\nabla \in \E \tensoran \Omega^2(A_N) \tensoran \E^*$.

In particular for the Levi-Civita connection of the fuzzy sphere with curvature in \eqref{curv} one easily finds,
$$
R_\nabla =
\tfrac{1}{4} \sum\nolimits_{p j} \, e_p \tensoran ( e_p \wedge e_j ) \, \tensoran \psi_j \, .
$$

Being $\E$ and $\E^*$ centered bimodules, we have a we well defined $ A_N - A_N $-bilinear map, 
\begin{align} \label{16thaugust2019n} 
\widehat{{\rm ev}} : \Omega^2(A_N) \tensoran \E^* \, \to \, \E \, & , \nn \\ 
 \widehat{{\rm ev}} \big( (e_j \wedge e_k) \tensoran \psi_m \big) 
 &=  \tfrac{1}{2} \, \widehat{{\rm ev}} 
 \big( (e_j \tensoran e_k - e_k \tensoran e_j) \tensoran \psi_m \big) \nn \\ 
& = e_j \, \psi_m(e_k) - e_k \, \psi_m(e_j) = e_j \, \delta_{m k} - e_k \, \delta_{m j}.
\end{align}

We are prepared to define the Ricci and scalar curvature of a connection.

		\bdfn
				The Ricci curvature $ {\rm Ric} $ is defined as the element in $ \E \tensoran \E $ given by
				\begin{equation} \label{23rdjan1} {\rm Ric}:= ( {\rm id}_\E \tensoran \widehat{{\rm ev}} ) ( R_\nabla) \in \E \tensoran \E \, . \end{equation}
Furthermore, if $ {\rm Ric} = \sum\nolimits_{k j} r_{kj} \, e_k \tensoran e_j$ for elements $ r_{kj } \in \mathcal{Z}(A_N)$, the scalar curvature {\rm Scal} is defined as:
				\begin{equation} \label{23rdjan2} 
				{\rm Scal}: = \sum\nolimits_{k j} r_{kj} \, g ( e_k \tensoran e_j ) \,  \in \, A_N. 
\end{equation}
		\edfn
			
\bppsn
The scalar curvature for the fuzzy sphere is $\frac{3}{4}$.  		
\eppsn		
			\noindent {\bf Proof:} By applying formula \eqref{16thaugust2019n} for the 2-forms $e_p \wedge e_j$ we first compute 
	$$
\sum\nolimits_{j} \widehat{{\rm ev}} ( (e_p \wedge e_j) \tensoran \psi_j ) = 3 e_p - e_p =   2 e_p . 
$$
When inserting in \eqref{curv}, this yields			
\begin{align*}
 {\rm Ric} &:= ( \id \tensoran \widehat{{\rm ev}}  ) [ \, \tfrac{1}{4} \sum\nolimits_{p j} \, e_p \tensoran ( e_p \wedge e_j ) \, \tensoran \psi_j \,  ] \\
&\:= \tfrac{1}{2} \, \sum\nolimits_{p} e_{p} \tensoran e_p . 
\end{align*}
Then, 
\begin{align*}
{\rm Scal} 
& = \tfrac{1}{2}\, \sum\nolimits_p \, g ( e_p \tensoran e_p ) \\	
& = \tfrac{3}{4} \, . 
\end{align*}

\vspace{4mm}

\noindent
{\bf Acknowledgement:} JB and DG were funded by a ``Research in Pairs" INDAM grant at ICTP, Trieste, and also by a ``Dipartimento di Eccellenza" fund of the Department of Mathematics and Geosciences, University of Trieste. DG was partially supported by JC Bose Fellowship and Grant from the Department of Science and Technology, Govt. of India. GL acknowledges partial support from INFN, Iniziativa Specifica GAST and from INDAM GNSAGA. Finally, JB and DG want to thank GL for the very kind hospitality during their stay in Trieste.

\vspace{1in}

\noindent
Jyotishman Bhowmick, Debashish Goswami: \\
Indian Statistical Institute, 
203, B. T. Road, Kolkata 700108, India\\
 jyotishmanb$@$gmail.com, goswamid$@$isical.ac.in\\

\noindent 
Giovanni Landi: \\ Matematica, Universit\`a di Trieste, Via A. Valerio, 12/1, 34127 Trieste, Italy \\
Institute for Geometry and Physics (IGAP) Trieste, Italy and INFN, Trieste, Italy \\
landi@units.it

\end{document}